\documentclass[preprint,12pt,authoryear]{elsarticle}

\usepackage{amsmath,amsfonts,amsthm,amssymb}
\usepackage{algorithmic}
\usepackage{algorithm}
\usepackage{array}
\usepackage[caption=false,font=small,labelfont=sf,textfont=sf]{subfig}
\usepackage{textcomp}
\usepackage{balance}
\usepackage{stfloats}
\usepackage{url}
\usepackage{verbatim}
\usepackage{graphicx}
\usepackage{color}

\usepackage{xspace}
\usepackage{bbm}
\usepackage{mathrsfs}

\def\diag{\mathop{\rm diag}\nolimits}%
\newcommand{\Hc}{\mathcal{H}}

\newcommand{\Rc}{\mathcal{R}}

\newcommand{\Wc}{\mathcal{W}}

\newcommand{\Av}{{\bf A}}
\newcommand{\Bv}{{\bf B}}
\newcommand{\Cv}{{\bf C}}

\newcommand{\Gv}{{\bf G}}

\newcommand{\Iv}{{\bf I}}

\newcommand{\Xv}{{\bf X}}

\newcommand{\Sv}{{\bf S}}

\newcommand{\bv}{{\bf b}}

\newcommand{\Ov}{{\bf O}}

\newcommand{\yv}{{\bf y}}
\newcommand{\wv}{{\bf w}}

\newcommand{\sv}{{\bf s}}

\newcommand{\betav}{\boldsymbol \beta}
\newcommand{\deltav}{\boldsymbol \delta}

\newcommand{\Sigmav}{\boldsymbol \Sigma}

\def\b{\beta}

\def\eps{\epsilon}
\def\l{\lambda}

\newcommand{\Norm}{\mathcal{N}}

\def\textiid{i.i.d.\@\xspace}
\newcommand\iid{\ifmmode\text{ i.i.d. } \else \textiid \fi}

\newcommand{\Real}{\mathbb{R}}

\newcommand{\ind}{\boldsymbol{1}}

\newcommand{\indep}{\perp \!\!\! \perp}

\newtheorem{theorem}{Theorem}
\newtheorem{lemma}{Lemma}

\newtheorem{assumption}{Assumption}

\newtheorem{corollary}{Corollary}
\newtheorem{remark}{Remark}

\usepackage{amssymb}

\allowdisplaybreaks
\begin{document}

\begin{frontmatter}



\title{Variable Selection with the Knockoffs: \\Composite Null Hypotheses}


\author{Mehrdad Pournaderi and Yu Xiang}
\affiliation{organization={Department of Electrical and Computer Engineering\\ University of Utah},
addressline=\\
            {50 Central Campus Dr 2110},
            city={Salt Lake City},
            postcode={84112}, 
            state={Utah},
            country={USA}}



\begin{abstract}
The fixed-X knockoff filter is a flexible framework for variable selection with false discovery rate (FDR) control in linear models with arbitrary design matrices (of full column rank) and it allows for finite-sample selective inference via the Lasso estimates. In this paper, we extend the theory of the knockoff procedure to tests with composite null hypotheses, which are usually more relevant to real-world problems. The main technical challenge lies in handling composite nulls in tandem with dependent features from arbitrary designs. We develop two methods for composite inference with the knockoffs, namely, shifted ordinary least-squares (S-OLS) and feature-response product perturbation (FRPP), building on new structural properties of test statistics under composite nulls. We also propose two heuristic variants of S-OLS method that outperform the celebrated Benjamini-Hochberg (BH) procedure for composite nulls, which serves as a heuristic baseline under dependent test statistics. Finally, we analyze the loss in FDR when the original knockoff procedure is naively applied on composite tests.\footnote{Publication link: https://doi.org/10.1016/j.jspi.2023.106119} 
\end{abstract}




\end{frontmatter}


\section{Introduction}

Selecting variables from a large collection of potential explanatory variables that are associated with responses of interest is a fundamental problem in many fields of science including genome-wide association study (GWAS), geophysics, and economics. In this paper, we focus on the classical linear regression model,
\begin{align}
	\yv = \Xv \betav + \wv\ ,\qquad \wv\sim \Norm(0, \sigma^2 \Iv_n)\ ,\label{linear_model}
\end{align}
where $\yv$ and $\wv$ are $n$-dimensional random vectors with elements denoting response and error variables, respectively, $\Xv = [\Xv_1,...,\Xv_p]\in\Real^{n\times p}$ denotes a fixed design matrix containing $n$ samples of $p$ explanatory features/variables, and $\betav = [\beta_1,...,\beta_p]^\top\in\Real^p$ is the vector of unknown fixed coefficients relating $\Xv$ and $\mathbb{E}(\yv)$. For the following $p$ hypotheses,
\begin{equation*}
\bigg\{
	\begin{array}{ll}
		\mathsf{H}_{0,i}: \beta_i=0\\
		\mathsf{H}_{1,i}: \beta_i\neq 0
	\end{array}
\ , \quad 1\leq i\leq p\,
\end{equation*}
the problem of interest is to test these hypotheses while controlling a simultaneous measure of type I error called the \emph{false discovery rate (FDR)} introduced by~\cite{benjamini1995}. {Let $\Hc_0=\{1\leq i\leq p:\beta_i=0\}$ denote the set of variables for which the null hypothesis is true and $\Rc=\{i:\mathsf{H}_{0,i} \text{ rejected}\}$ denote the selected variables by some variable selection procedure. The FDR is defined as} $\mathsf{FDR}=\mathbb{E}(\mathsf{FDP})$ with
\begin{equation}
     \mathsf{FDP}=\frac{|\Rc\cap\Hc_0|}{|\Rc|\vee 1},\label{def:fdr}
\end{equation}
where $|\cdot|$ denotes the cardinality of the sets. A selection rule controls the FDR at level $q$ if its corresponding FDR is guaranteed to be at most $q$ for some predetermined $q\in [0,1]$. 

Recently, \cite{barber2015controlling} proposed the (fixed-X) knockoff filter procedure, a data-dependent selection rule that controls the FDR in finite sample settings and under arbitrary designs. In this procedure, a test statistic is computed for each feature through constructing a knockoff variable, and a feature is selected by (data-dependent) thresholding the statistics according to the target FDR. The knockoff construction allows for correlated features {and \cite{barber2015controlling} considers a range of simulations settings for which} the knockoff framework has higher statistical power in comparison with the Benjamini-Hochberg (BH) procedure introduced in \cite{benjamini1995} (see \cite{benjamini2001control,storey2004strong,efron2001empirical} for extensions and variants).
The knockoff filter has inspired various formulations such as the model-X knockoffs and deep learning-based knockoffs among others (see \cite{candes2018panning,barber2019knockoff,barber2020robust,romano2019deep,jordon2018knockoffgan,lu2018deeppink,fan2019ipad,pournaderi2021differentially}). 

The fixed-X knockoff filter was originally designed to test simple null hypotheses ($\beta_i=0$), but in practice (e.g. microarray experiment and GWAS), researchers often encounter composite null hypotheses. Composite null hypotheses are relevant for two reasons. Firstly, the traditional sparsity assumption in linear models, which assumes that only a few coefficients are non-zero, may not be valid. Instead, a more realistic null hypothesis for a variable is that it has little or negligible effect on the response variable. Secondly, 
researchers are often interested in detecting ``large" effects, and non-zero coefficients with a relatively weak effect may not be of interest for detection. Both of these issues have been discussed in detail in the previous literature and we refer to \cite{barber2019knockoff,doi:10.1080/01621459.2012.664505} and references therein for further discussion.
The multiple testing of composite null hypotheses using independent p-values has been studied in \cite{benjamini2001control,doi:10.1080/01621459.2012.664505,DICKHAUS20131968,CABRAS2010659}. 

Given that the knockoff selection framework deals with the dependencies between statistics inherently, a natural question is whether one can extend this to handle composite nulls, namely, 
\begin{equation*}
\bigg\{
	\begin{array}{ll}
		\Breve{\mathsf{H}}_{0,i}: |\beta_i|\leq {\delta_i}\\
		\Breve{\mathsf{H}}_{1,i}:  |\beta_i|>{\delta_i}
	\end{array}
\ , \quad 1\leq i\leq p\ ,\label{comptest}
\end{equation*}
for some given {set of $\delta_i \geq 0$}.
In this paper, we provide an affirmative answer to the above question by developing two methods: shifted ordinary least-squares (S-OLS) and feature-response product perturbation (FRPP). We show that both methods achieve FDR control in finite sample settings under arbitrary designs, leveraging new structural properties for test statistics under composite nulls. {Furthermore, we compute an upper bound for the FDR when one uses the original fixed-X knockoff filter for composite nulls, and the result reduces to the exact FDR control for simple nulls.} The main technical difficulty is to handle composite nulls in tandem with \emph{dependent} features. 
The original BH procedure can handle composite nulls under the independence
assumption on the test statistics~\cite[Theorem 5.2]{benjamini2001control}.
In order to handle arbitrarily dependent test statistics, \cite{benjamini2001control} proposed a quite conservative correction for the test size, and the method is known as the BY procedure. The BY procedure has later been shown in~\cite{blanchard2008two} to be theoretically valid for composite tests (super-uniform p-values) as well. Another existing solution to the problem is to apply the BH procedure on the p-values obtained from a knockoff-assisted estimation of $\betav$ which provides independent estimates for the model coefficients~\cite{barber2019knockoff,sarkar2022adjusting}. We refer to this method as the knockoff-assisted BH procedure. The aforementioned approach, however, results in a substantial power loss in comparison to the direct (heuristic) use of the BH procedure.
In a recent study, \cite{fithian2022conditional} has attempted to address the dependent test statistics without sacrificing the statistical power of the BH procedure, and proposed the dependence-adjusted Benjamini–Hochberg (dBH) procedure. The study shows empirically that the dBH method performs similarly to BH in terms of power, but with provable FDR control. 
This method can handle one-sided composite nulls in exponential family models and the paper claims that the method is extensible to the two-sided composite null hypotheses as well. The conditional calibration framework of~\cite{fithian2022conditional} has been adopted in~\cite{luo2022improving} to improve the knockoff procedure.
 
{In the simulations section, it is shown that the FRPP knockoff procedure, a randomized generalization of the knockoff procedure, outperforms the knockoff-assisted BH procedure which serves as the knockoff-based theoretical baseline in this paper. Our S-OLS method exhibits similar performance as the BY procedure. Also, we use the BH procedure (without any correction) 
as a heuristic baseline due to its wide use in applications.} The S-OLS method motivates two Lasso-based heuristic variants that outperform the composite BH procedure in power.  

The paper is organized as follows. In Section~\ref{sec:background}, we briefly present the knockoff filter framework {of~\cite{barber2015controlling}} by introducing the main steps to set the stage for our analysis. In Section~\ref{sec:main}, we present our main results and theoretical guarantees, along with two heuristic methods, with all the proofs deferred to Appendices. We report our experimental results in Section~\ref{sec:sim} for a range of composite nulls, amplitude of alternatives, and correlation coefficients. 

\section{Background: Fixed-X Knockoff Filter}
\label{sec:background}

The knockoff variable selection procedure \cite{barber2015controlling} consists of two main steps: (I) computing a statistic $W_j$ for each variable in the model, and (II) selecting ${\mathcal{R}}_{ko}=\{j:W_j\geq T(q,\Wc)\}$, where the threshold $T\geq 0$ depends on the target FDR $q$ and the set of computed statistics, i.e., $\Wc=\{W_1,\hdots,W_p\}$. 
In this section, we look at this procedure in more detail.
\subsection{Knockoff Design}
The knockoff methodology for detecting non-null variables involves creating a {control} (or knockoff) design that mimics the correlation structure of $\Xv$ but {its relationship to $\yv$ does not (necessarily) follow the one between $\Xv$ and $\yv$.} 
 \begin{assumption}
 $\Sigmav=\Xv^\top\Xv$ is invertible.
 \end{assumption}
Specifically, if $n\geq 2p$, \cite{barber2015controlling} proposes the following construction to produce knockoff designs   
\begin{equation}
    \mathbf{\Tilde{X}}(\sv)=\Xv\,({\bf I}_p-\boldsymbol{\Sigma}^{-1} \text{diag}\{\sv\})+\mathbf{\tilde{U}}{\bf C},\label{KOconst}
\end{equation}
where $\sv\in\Real^p_+$ is a free vector of parameters as long as it satisfies {$\diag(\sv)\preceq 2\Sigmav$ which guarantees the existence of the Cholesky decomposition $\Cv^\top\Cv=2\diag(\sv)-\diag(\sv)\Sigmav^{-1}\diag(\sv)$;}
 and $\mathbf{\tilde{U}}^{n\times p}$ is an orthonormal matrix that satisfies $\mathbf{\tilde{U}}^\top\Xv=0$ (see \cite{barber2015controlling} for details). Therefore,
knockoff matrices are not unique and they are constructed according to the original design matrix. Using \eqref{KOconst}, the following relation can be easily verified.
\begin{align}
    {\Gv := [ \Xv\ \mathbf{\Tilde{X}}]^\top[ \Xv\ \mathbf{\Tilde{X}}]} = 
    \begin{bmatrix} \boldsymbol{\Sigmav} &  \boldsymbol{\Sigmav}-\diag(\sv)\\\boldsymbol{\Sigmav}-\diag(\sv) &  \boldsymbol{\Sigmav}\end{bmatrix}.\label{eq:data_gram}
\end{align}
In fact, this construction not only preserves the correlation structure of $\Xv$, but also has another subtle yet important geometrical implication: $\Xv_i$ ($i$-th column of $\Xv)$ and $\mathbf{\Tilde{X}}_i$ are second-order ``exchangable'' in a deterministic sense, i.e., swapping $\Xv_i$ and $\mathbf{\Tilde{X}}_i$ does not change the inner product structure (Gram matrix) of the augmented design. This property makes $\mathbf{\Tilde{X}}$ an appropriate tool for FDR control.
It should be noted that the detection power highly depends on the parameter $\sv$ as it determines the angle between a feature and its corresponding knockoff. In other words, $s_i$ ($i$-th element of $\sv$) controls how different (or orthogonal) $\Xv_i$ and $\mathbf{\Tilde{X}}_i$ would be. Assuming the columns of $\Xv$ are normalized by the Euclidean norm, one way to choose $\sv$ is to solve the following convex problem,
\begin{align}
     \text{minimize} \qquad &\sum_{i=1}^p |1-s_i|\label{sdpknockoffs}\\
     \text{subject to} \qquad & s_i \geq 0\ ,\nonumber
     \quad\text{diag} \{\sv\} \preceq 2\,\Sigma \ .
\end{align}
This semi-definite programming minimizes the average correlation between variables and their corresponding knockoff variable. {See~\cite{barber2015controlling,spector2022powerful} for other choices of $\sv$}.
\begin{remark}
Since a column-wise normalization of the design matrix is natural for variable selection purposes and essential in terms of statistical power, throughout this paper, we always assume that $\Xv$ is normalized by the $\ell_2$ norm of columns.
\end{remark}
\subsection{Statistics}
Using the knockoff features, we now compute a vector of anti-symmetric statistics $\boldsymbol{W}=(W_1, W_2,\hdots,W_p)^\top$ by regression over the augmented design $[\Xv\ \mathbf{\Tilde{X}}]$. Let $\hat\theta_i$ and $\hat\theta'_i$ denote some estimated parameters corresponding to the variables $X_i$ and $\tilde{X}_i$. In this case, one can define a statistic as follows
\begin{equation}
     W_i = |\hat\theta_i| - |\hat\theta'_i|\ ,\qquad 1\leq i\leq p\ .\label{stat:LCD}
\end{equation}
We can also define the statistics differently,
\begin{equation}
    W_i = \mathsf{sgn}(|\hat\theta_i| - |\hat\theta'_i|)\,\max(|\hat\theta_i|,|\hat\theta'_i|)\ .\label{lcsm}
\end{equation} To be more precise, the term \textit{anti-symmetric} here means that swapping the estimates $\hat\theta_i$ and $\hat\theta'_i$ for any subset of indices $F\subseteq \{1,2,\hdots,p\}$ has the effect of switching the signs of $\{W_i:i\in F\}$. Specifically, the knockoff framework guarantees the FDR control when the (anti-symmetric) statistics are computed based on estimators that depend on the data through the following form,
\begin{equation}
    \begin{pmatrix}\boldsymbol{\hat\theta}_{p\times 1}\\\boldsymbol{\hat\theta'}_{p\times 1}\end{pmatrix}= \mathcal{E}(\Gv,[\Xv\ \mathbf{\Tilde{X}}]^\top\yv)\label{general_form} \ ,
\end{equation}
{where $\mathcal{E}$} is a deterministic operator and, swapping $\Xv_i$ and $\mathbf{\Tilde{X}}_i$ will result in swapping $\hat\theta_i$ and $\hat\theta'_i$.
For instance, the Lasso~\cite{tibshirani1996regression} regression estimates $\boldsymbol{\hat{\beta}}_\mathsf{LASSO}$ given by
\begin{equation}
    \boldsymbol{\hat{\beta}}_\mathsf{LASSO}:=\underset{\bv\,\in\Real^{2p}}{\arg\min}\Big\{\,\Big\|\yv-[\Xv\ \mathbf{\Tilde{X}}]\bv\Big\|_2^2\,+\lambda \|\bv\|_1 \Big\},\quad \lambda\geq 0,\label{eq:LS_model}
\end{equation}
can be considered as an example of estimators satisfying $\eqref{general_form}$.
\subsection{FDR Control}
In this subsection, we briefly discuss the existing approaches to show the FDR control in the knockoff procedure.  Let $P_F$ denote the $2p\times 2p$ permutation matrix corresponding to swapping $\Xv_i$ and $\mathbf{\Tilde{X}}_i$ for all $i\in F$, then by the structure of knockoff matrix we have,
\begin{equation}
    P_F^\top\Gv P_F = \Gv\ ,\qquad F\subseteq\{1,\hdots,p\}\ .\label{exch1}
\end{equation}
Also, according to \eqref{linear_model} we get
\begin{equation}
    P_F^\top[ \Xv\ \mathbf{\Tilde{X}}]^\top\yv \overset{d}{=} [\Xv\ \mathbf{\Tilde{X}}]^\top\yv\ ,\qquad F\subseteq\Hc_0\ . \label{exch2}
\end{equation}
The identities \eqref{exch1} and \eqref{exch2} immediately imply an interesting property of the estimates $\boldsymbol{\hat{\theta}}$ and $\boldsymbol{\hat{\theta}'}$: \textit{the estimated parameters for null variables and their corresponding knockoff variables are {(conditionally)} exchangeable.} This property along with the anti-symmetric structure of the statistics are the main ingredients of the FDR control proof in the original paper.   
In fact, under the simple null hypotheses, these properties lead to {conditionally symmetric statistics, i.e., $W_j|W_{-j}\overset{d}{=}-W_j|W_{-j}$ for all $j\in\Hc_0$. This symmetry property is called \textit{\iid sign property} for the nulls in~\cite{barber2015controlling}, meaning that} the signs of the null statistics (signs of $\Wc_0=\{W_j: \beta_j=0\}$) are independent of the magnitudes and have \iid Rademacher distribution (i.e., $\mathbb{P}(+1)=\mathbb{P}(-1)=1/2$). In this case, using the martingale theory, it is shown that rejecting $\{j:W_j\geq T\}$ with the following threshold $T$ controls the FDR at level $q$,
\begin{equation}
    T = \inf\Big\{t \in \Phi:\widehat{\mathsf{FDP}}(t)\leq q\Big\},\label{threshold}
\end{equation}
\begin{equation}
    \widehat{\mathsf{FDP}}(t):=\frac{1+\#\{j: W_j \leq -t\}}{\#\{j: W_j\geq t\}\vee 1},
\end{equation}
where $\Phi =\{|W_j|:j=1,2,\hdots p\}\setminus\{0\}$.
Although this approach is elegant, we shall not follow it in the development of composite tests as in this situation \eqref{exch2} fails to hold immediately. On the other hand, \cite{barber2020robust} provides another proof for FDR control\footnote{The results in \cite{barber2020robust} concern robustness of the Model-X knockoff framework (where $\Xv$ is random) but the FDR control proof works for the fixed-X setting as well since they only rely on antisymmetry of the statistics and \eqref{robustenssineq}.} which requires weaker conditions on the statistics \cite[Equation~(16)]{barber2020robust}. Specifically, it suggests that if the statistics {corresponding to true null hypotheses} satisfy 
\begin{align}
    &\mathbb{P}\Big\{W_j > 0 \ \Big|\ |W_j|, \boldsymbol{W}_{-j}\Big\}\leq C\cdot \mathbb{P}\Big\{W_j < 0 \ \Big|\ |W_j|, \boldsymbol{W}_{-j}\Big\}\label{robustenssineq},
\end{align}
{(almost surely)} for some $C>0$, then the knockoff procedure with the target FDR $q/C$ controls the FDR at level $q$. It is straightforward to verify that in the case of simple nulls, this condition holds with $C\geq 1$ according to the \iid sign property for the nulls, resulting in FDR control.

\section{Main Results}
\label{sec:main}
In the case of a single composite test of the form $|\beta_i|\leq\delta_i$, the common approach is to compute super-uniform p-values under the null (i.e., $\mathbb{P}\,(P\leq t)\leq t$) which clearly controls the probability of type I error by the definition. The super-uniformity usually happens when a p-value is computed according to the distribution corresponding to a parameter on the boundary of the null region. However, for composite multiple testing problems {under the FDR control constraint}, this argument gets more complicated as the dependencies between the statistics should be considered. {The BH procedure guarantees the FDR control when the null p-values are super-uniform, mutually independent, and independent of the non-null p-values according to \cite{ramdas2019unified}.} On the other hand, the knockoff filter is designed to utilize the model and covariates structure for computing one-bit p-values that handle the dependencies naturally. It turns out that one can actually maintain this property of the knockoff procedure when developing composite tests. {Let $\Breve\Hc_0 = \{1\leq i\leq p:|\beta_i|\leq \delta_i\}$ denote the set of indices for which the composite null hypothesis is true.} To show the FDR control, we rely on manipulating the estimates so that the {null statistics ($j\in\Breve\Hc_0$)} satisfy the following inequality for some $B<\infty$,
\begin{align*}
    \mathbb{P}\Big\{W_j > 0 \ \Big|\ |W_j|, \boldsymbol{W}_{-j}\Big\}\leq B\cdot \mathbb{P}\Big\{W_j < 0 \ \Big|\ |W_j|, \boldsymbol{W}_{-j}\Big\}\ a.s.,\ \text{all } j\in \Breve\Hc_0,
\end{align*}
which is equivalent to
\begin{equation}
    R:={\underset{j \in \Breve\Hc_0}{\max}\ \underset{\omega:W_j \neq 0}{\text{ess}\sup}}\ \frac{\mathbb{P}\Big\{W_j > 0 \ \Big|\ |W_j|, \boldsymbol{W}_{-j}\Big\}}{\mathbb{P}\Big\{W_j < 0 \ \Big|\ |W_j|, \boldsymbol{W}_{-j}\Big\}}\leq B . \label{FDRGbound}
\end{equation}
It is known from \cite[Theorem~2 with $E_j=0$]{barber2020robust} that having this bound leads to rigorous FDR control at level $q\cdot B$. In fact, showing $R\leq 1$ in the knockoff procedure framework means that the procedure overesitmates ${\mathsf{FDP}}(t)$ and therefore, can be interpreted as an equivalent for super-uniformity of null p-values in terms of BH procedure. However, for bounds $R\leq B$ where $B>1$, one needs to correct the test size by a factor of $1/B$, i.e., $q'=q/B$.
In cases where $R=\infty$, we use a generalization of this argument which leads to tighter bounds. We bound the following quantity.
\begin{equation}
    R':={\underset{j \in \Breve\Hc_0}{\max}\ \underset{\omega:W_j \neq 0}{\text{ess}\sup}}\ \frac{\mathbb{P}\Big\{W_j > 0,\, A_j \ \Big|\ |W_j|, \boldsymbol{W}_{-j}\Big\}}{\mathbb{P}\Big\{W_j < 0 \ \Big|\ |W_j|, \boldsymbol{W}_{-j}\Big\}}, \label{FDRGbound2}
\end{equation}
where $A_j$ is some event regarding the $j$-th variable. 
\begin{theorem}[\cite{barber2020robust}, Theorem 2]
\label{bcthm}
If $R'\leq B'$, we get $\mathsf{FDR}\leq q\cdot B'+\mathbb{P}\big(\bigcup_{j\in\Breve\Hc_0}A_j^c\big)$.
\end{theorem}
In the following section, we present composite selective inference methods that allow for theoretical FDR control. 
\subsection{Composite Testing with FDR Control}
The following theorem concerns the composite knockoff procedure based on the ordinary least-squares estimates\footnote{Note that performing the knockoff procedure using the OLS estimator requires that $\Gv$ is invertible. In Lemma \ref{invlem} it is shown that this is the case if $\underset{1\leq i\leq p}{\max}\,s_i<2\lambda_{\min}(\Sigmav)$.} {$\boldsymbol{\hat{\beta}}_\mathsf{OLS}=\begin{pmatrix}\boldsymbol{\hat\beta}\\\boldsymbol{\hat\beta'}\end{pmatrix}=\Gv^{-1}[\Xv\ \mathbf{\Tilde{X}}]^\top\yv$}. We consider both one-sided and two-sided null hypotheses, i.e., $\beta_j\leq \delta_j$ and $|\beta_j|\leq \delta_j$, respectively. We show that shifting the estimates corresponding to the knockoff variables by $\delta_j$ will result in exact FDR control for one-sided test {(i.e., we prove $R\leq 1$ in this case) and extend it to a two-sided test via a Bonferroni correction argument. We also propose a method with approximate FDR control for two-sided tests. In this case}, we derive a bound for the FDR of the knockoff procedure {(i.e., $R=\infty$ and we show $R'\leq 1$)}.
\begin{remark}
    Let $D=\diag(\sv)$. It is noteworthy that 
    \[
    \boldsymbol{\hat\beta}-\boldsymbol{\hat\beta'}=D^{-1}(\Xv-\mathbf{\Tilde{X}})^\top\yv\sim \Norm(\betav,2\sigma^2 D^{-1}),
    \]
    coincide with one of the two estimators investigated in \cite{sarkar2022adjusting}. 
\end{remark}

\begin{theorem}[\textbf{S-OLS}]
\label{thm:SOLS}
Consider the knockoff procedure with target FDR $q$ and based on the estimates $\boldsymbol{\hat{\beta}}_\mathsf{S-OLS}=\boldsymbol{\hat{\beta}}_\mathsf{OLS}+\begin{pmatrix}\mathbf{0}_{p\times 1}\\\deltav_{p\times1}\end{pmatrix}$. Let $\hat\beta_j$ and $\hat\beta_j'$ denote the $j$-th and $(j+p)$-th elements of $\boldsymbol{\hat{\beta}}_\mathsf{S-OLS}$.
\medskip

\noindent
\textbf{(I) One-sided test:} Consider testing {$\mathsf{H}'_{0,j}:\beta_j\leq \delta_j$ for $1\leq j\leq p$. If $|W_j|$ depends on the estimator $\boldsymbol{\hat{\beta}}_\mathsf{S-OLS}$ through the unordered pair $\{{\hat\beta}_j,{\hat\beta}_j'\}$ and $\mathsf{sgn}(W_j) = \mathsf{sgn}\big({\hat\beta}_j-{\hat\beta}_j'\big)$, then the procedure controls the FDR at level $q$.}
\smallskip

\noindent
\textbf{(II) Approximate two-sided test:} Consider testing {$\Breve{\mathsf{H}}_{0,j}:|\beta_j|\leq \delta_j$ for $1\leq j\leq p$. If $|W_j|$ depends on the estimator $\boldsymbol{\hat{\beta}}_\mathsf{S-OLS}$ through the unordered pair $\{{\hat\beta}_j,{\hat\beta}_j'\}$} and $\mathsf{sgn}(W_j) = \mathsf{sgn}\big(\big|{\hat\beta}_j\big|-\big|{\hat\beta}_j'\big|\big)$, then 
\begin{align*}
    \mathsf{FDR} \leq & \, q + \mathbb{P}\bigg\{\underset{j\in\Breve\Hc_0}{\min}\,\big(\hat\beta_j+\hat\beta_j'\big)<0\bigg\}
    \,.
\end{align*}\label{shifted_OLS}
\end{theorem}
{For an exact two-sided test using the OLS estimator, one can consider testing the intersection hypothesis $\Breve{\mathsf{H}}_{0,j}:\mathsf{H}'_{0,j}\cap \mathsf{H}''_{0,j}$, $1\leq j\leq p$, via a Bonferroni-type method, where $\mathsf{H}'_{0,j}:\beta_j\leq \delta_j$ and $\mathsf{H}''_{0,j}:\beta_j\geq -\delta_j$. To be precise, one needs to perform the knockoff procedure with target FDR $q/2$ twice; once to test $\mathsf{H}'_{0,j}:\beta_j\leq \delta_j$ and another one for $\mathsf{H}''_{0,j}:\beta_j\geq -\delta_j$. Let $\Rc'_{q/2}=\{j:\mathsf{H}'_{0,j} \text{ rejected}\}$ and define $\Rc''_{q/2}$ similarly. The two-sided procedure rejects the union $\Rc=\Rc'_{q/2}\cup \Rc''_{q/2}$. In order to test $\mathsf{H}''_{0,j}:\beta_j\geq -\delta_j$, the same procedure as Theorem \ref{shifted_OLS} (I) can be used but with $\boldsymbol{\hat{\beta}}_\mathsf{S-OLS}=-\boldsymbol{\hat{\beta}}_\mathsf{OLS}+\begin{pmatrix}\mathbf{0}\\\deltav\end{pmatrix}$. Let $\Hc_0'=\{j:\mathsf{H}'_{0,j} \text{ true}\}$ and define $\Hc_0''$ similarly.
\begin{corollary}[\textbf{Exact two-sided test}]
Consider testing the two-sided null hypotheses $\Breve{\mathsf{H}}_{0,j}:|\beta_j|\leq \delta_j$ for $1\leq j\leq p$. Rejecting $\Rc=\Rc'_{q/2}\cup \Rc''_{q/2}$ controls the FDR at level $q$, i.e.,
\begin{equation}
     \mathsf{FDR} =\mathbb{E}\left(\frac{|\Rc\cap\Breve{\Hc}_0|}{|\Rc|\vee 1}\right)\leq q,
\end{equation}   
where $\Breve{\Hc}_0=\Hc_0'\cap \Hc_0''$.
\begin{proof}
    \begin{align*}
         \mathbb{E}\left(\frac{|\Rc\cap\Breve{\Hc}_0|}{|\Rc|\vee 1}\right)&= \mathbb{E}\left(\frac{|(\Rc'_{q/2}\cap\Breve{\Hc}_0)\cup(\Rc''_{q/2}\cap\Breve{\Hc}_0)|}{|\Rc|\vee 1}\right)\\
         &\leq \mathbb{E}\left(\frac{|\Rc'_{q/2}\cap\Breve{\Hc}_0|}{|\Rc|\vee 1}\right)+\mathbb{E}\left(\frac{|\Rc''_{q/2}\cap\Breve{\Hc}_0|}{|\Rc|\vee 1}\right)\\
         & \leq \mathbb{E}\left(\frac{|\Rc'_{q/2}\cap \Hc'_0|}{|\Rc'_{q/2}|\vee 1}\right)+\mathbb{E}\left(\frac{|\Rc''_{q/2}\cap\Hc''_0|}{|\Rc''_{q/2}|\vee 1}\right) \leq \frac{q}{2} + \frac{q}{2} = q.
    \end{align*}
\end{proof}
\end{corollary}
}
\begin{corollary}[\textbf{Alternative method for Theorem \ref{shifted_OLS} (II)}] \label{corol_abs} 
Under the same setting as Theorem \ref{shifted_OLS} (II) we get,
\begin{align*}
    \mathsf{FDR} \leq & \, q + \mathbb{P}\bigg\{\underset{j\in\Breve\Hc_0}{\max}\,\big(\hat\beta_j+\hat\beta_j'\big)>0\bigg\}
    \, ,
\end{align*}
{when $\boldsymbol{\hat{\beta}}_\mathsf{S-OLS}=\boldsymbol{\hat{\beta}}_\mathsf{OLS}-\begin{pmatrix}\mathbf{0}_{p\times 1}\\\deltav_{p\times1}\end{pmatrix}$.}
\end{corollary}

The next method generalizes the knockoff framework in the sense that composite inference is allowed but it is not limited to any particular estimator. In this situation, the shifted estimates are no longer feasible to analyze. Therefore, in order to perform composite inference in such a general setting, we propose to introduce artificial randomness to the procedure. In particular, we perturb the feature-response products $[\Xv\ \mathbf{\Tilde{X}}]^\top\yv$ by noise generated from Laplace distribution. In this case, we will be able to show that $R\leq e^\eps$ with $\eps$ being determined by the noise variance.

\begin{theorem}[\textbf{FRPP}]
\label{thm:FRPP}
Fix some $\eps >0$. Define the null variables $\Breve\Hc_0 = \{i:|\beta_i|\leq \delta_i\}$ and let {$\Delta_j{\sim}\text{Lap}\,(2\,s_{j\,(\mathsf{mod}\,p)}\,\delta_{j\,(\mathsf{mod}\,p)}/\eps)$, $1\leq j\leq 2p$, be $2p$ independent random variables} where $s_i = 1-\Xv_i^\top\mathbf{\Tilde{X}}_i$, $1\leq i\leq p$. The knockoff procedure with the target FDR $q/e^\eps$, and using (antisymmetric) statistics based on any estimator of the form $\boldsymbol{\tilde{\theta}}=\mathcal{E}(\Gv,[\Xv\ \mathbf{\Tilde{X}}]^\top\yv+\Delta)$ controls the FDR at level $q$. \label{lap_mech}
\end{theorem}
\begin{remark}
As we discussed in the previous section, the original knockoff framework focuses on the estimators of the form \eqref{general_form} that satisfy $P_i\,\boldsymbol{\hat{\theta}} = \mathcal{E}(\Gv,P_i[\Xv\ \mathbf{\Tilde{X}}]^\top\yv)$ where $P_i$ is the (symmetric) permutation matrix swapping the $i$-th and $(i+p)$-th elements. We also keep assuming this property as we refer to general estimators, e.g., we assume $P_i\,\boldsymbol{\tilde{\theta}}=\mathcal{E}\big(\Gv,P_i([\Xv\ \mathbf{\Tilde{X}}]^\top\yv+\Delta)\big)$. Observe that this is a very mild assumption (since $P_i\Gv P_i = \Gv$) and is satisfied by almost every estimator of linear models.
\end{remark}
We note that Theorem \ref{lap_mech} gives a stochastic generalization of the knockoff procedure, i.e., if {$\delta_i = 0$ for all $1\leq i \leq p$, then} the FRPP knockoff procedure reduces to the original method without any additional assumptions. 

\subsection{Heuristic Methods}
Motivated by our results that show the shifting argument is theoretically valid in the case of using the OLS estimator, we propose the following two methods based on shifting the Lasso estimates.

\noindent\textbf{Method {S-LASSO1}:} This method shifts the Lasso estimates just as the S-OLS method. Namely, we use the following formulae to compute the statistics.
\begin{equation*}
    \boldsymbol{\hat{\beta}}_\mathsf{S-LASSO1}=\boldsymbol{\hat{\beta}}_\mathsf{LASSO}(\lambda)+\begin{pmatrix}0\\\deltav\end{pmatrix}
\end{equation*}

\noindent\textbf{Method {S-LASSO2}:} This method estimates the coefficients by solving the following Lasso problem.
\begin{equation*}
    \boldsymbol{\hat{\beta}}_\mathsf{S-LASSO2}:=\underset{\bv\,\in\Real^{2p}}{\arg\min}\Bigg\{\,\bigg\|\yv-\big[\Xv\ \mathbf{\Tilde{X}}\big]\Big(\bv - \begin{pmatrix}0\\\deltav\end{pmatrix}\Big)\bigg\|_2^2\,+\lambda \big\|\bv\big\|_1 \Bigg\},\label{slas2}
\end{equation*}
where $\lambda\geq 0$. 
\begin{remark}
We note that both methods reduce to the S-OLS method if $\lambda = 0$.
\end{remark}

\subsection{FDR Bound for Naive Selection}
\begin{theorem}\label{boundthm}
Naive application of the fixed-X knockoff procedure with target FDR $q$ on composite null hypotheses $\Breve\Hc_0 = \{i:|\beta_i|\leq \delta_i\}$ and using (antisymmetric) statistics based on any estimator of the form $\mathcal{E}(\Gv,[\Xv\ \mathbf{\Tilde{X}}]^\top\yv)$, will result in FDR  bounded as follows.
\begin{equation}
    \mathsf{FDR} \leq \underset{\eps\geq 0}{\min}\Bigg\{q\,.\,e^\eps + \mathbb{P}\bigg({\frac{1}{\sigma^2}\,\underset{j\in\Breve\Hc_0}{\max}\,\left\{\delta_j\big|\gamma_j-\gamma_j'\big|\right\}}>\eps\bigg)\Bigg\}\ , \label{naivebnd}
\end{equation}
where $\begin{pmatrix}\boldsymbol{\gamma}_{p\times 1}\\\boldsymbol{\gamma'}_{p\times 1}\end{pmatrix}=[\Xv\ \mathbf{\Tilde{X}}]^\top\yv$.
\end{theorem}

\begin{remark}
We note that the bound \eqref{naivebnd} will reduce to $q$ in the case of simple nulls, i.e., {$\delta_i = 0$ for all $1\leq i \leq p$}.
\end{remark}

\section{Simulations}
\label{sec:sim}
 In this section, we present simulation results on synthetic data sets for all methods. We set the sample size and dimension to be $n=2000$ and $p = 800$, respectively. The samples (rows of $\Xv$) are generated \iid according to $\Norm(0,\Sv_p)$ where $\Sv_{ij}=\rho^{|i-j|},\ |\rho|<1$ and we normalize the columns of $\Xv$ by the $\ell_2$-norm. The responses $\yv$ are generated according to the linear model \eqref{linear_model} with noise variance $\sigma^2 = 1$ and the number of false nulls is $k=100$. The composite null boundary is set to be {$\delta_i=\delta=1$ for all $1\leq i \leq p$}, and we consider two distributions for generating coefficients corresponding to null variables: 
 
 \textbf{(a)} $\beta_0\sim U[-1,1]$ which is presented in the left column of the figures. We consider this setting as a practical case.
 
 \textbf{(b)} $\beta_0\sim \text{Rademacher}$, which is presented in the right column of the figures. This setting tries to examine the methods in the hardest situation (worst-case scenario) for FDR control.
 
We adopt the equicorrelated knockoffs and set the elements of $\sv$ (the vector used in creating the knockoff matrices; see \eqref{KOconst}) to be $\min(1.8\lambda_{\min}(\Sigmav),1)$, $\min(\lambda_{\min}(\Sigmav),1)$, and $\min(2\lambda_{\min}(\Sigmav),1)$ for the S-OLS, FRPP, and heuristic methods (S-LASSO1 and S-LASSO2), respectively. {More sophisticated constructions are available, including SDP knockoffs~\eqref{sdpknockoffs} and MRC knockoffs of \cite{spector2022powerful}, which can potentially improve the power of the knockoff-based methods in our simulations.} We adopt the structure \eqref{lcsm} to construct the coefficient signed max statistics and use the Lasso estimator \eqref{eq:LS_model} {to perform the FRPP knockoff procedure. The Lasso estimator is always used} with shrinkage parameter $\lambda=1$. {We do not have yet a valid cross-validation procedure under composite nulls for choosing $\lambda$ in the Lasso estimator. Therefore, the choice of $\lambda$ should be prespecified and we have chosen $\lambda = \delta = 1$, heuristically.}  

{We consider the BY procedure (which is the same as the BH procedure but with corrected test size $q/(\sum_{i=1}^p {1/i})$) via $\boldsymbol{\hat{\beta}}_\mathsf{BY}=\Sigmav^{-1}\Xv^\top\yv$ and the knockoff-assisted BH procedure via $\boldsymbol{\hat{\beta}}_\mathsf{KA}=\diag(\sv)^{-1}(\Xv-\mathbf{\Tilde{X}})^\top\yv$ (from~\cite{sarkar2022adjusting}) as theoretical baselines where two-sided p-values are computed according to 
\[
    P_{BY,j} = 2\Phi\left((\delta-|\hat\beta_{BY,j}|)/\sqrt{\sigma^2(\Sigmav^{-1})_{jj}}\right)
\]
and
\[
    P_{KA,j} = 2\Phi\left((\delta-|\hat\beta_{BY,j}|)\cdot\sqrt{s_j/(2\sigma^2)}\right)
\]
 with $\Phi(\cdot)$ denoting the CDF of the standard normal distribution and $\sigma^2=1$ is assumed to be known.}
 
\begin{figure}[h!]
\centering
\subfloat[][]{
\includegraphics[width=0.4\textwidth]{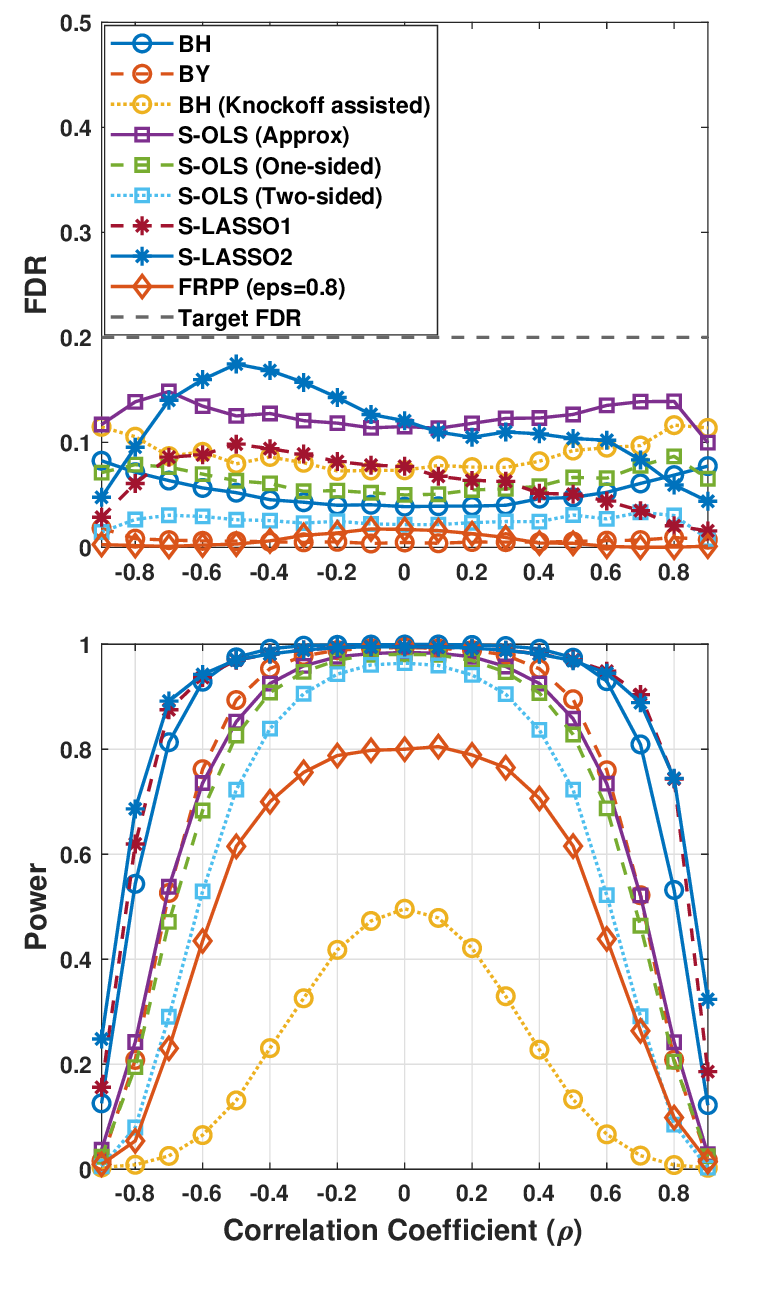}
\label{fig:subfig1}}
\subfloat[][]{
\includegraphics[width=0.4\textwidth]{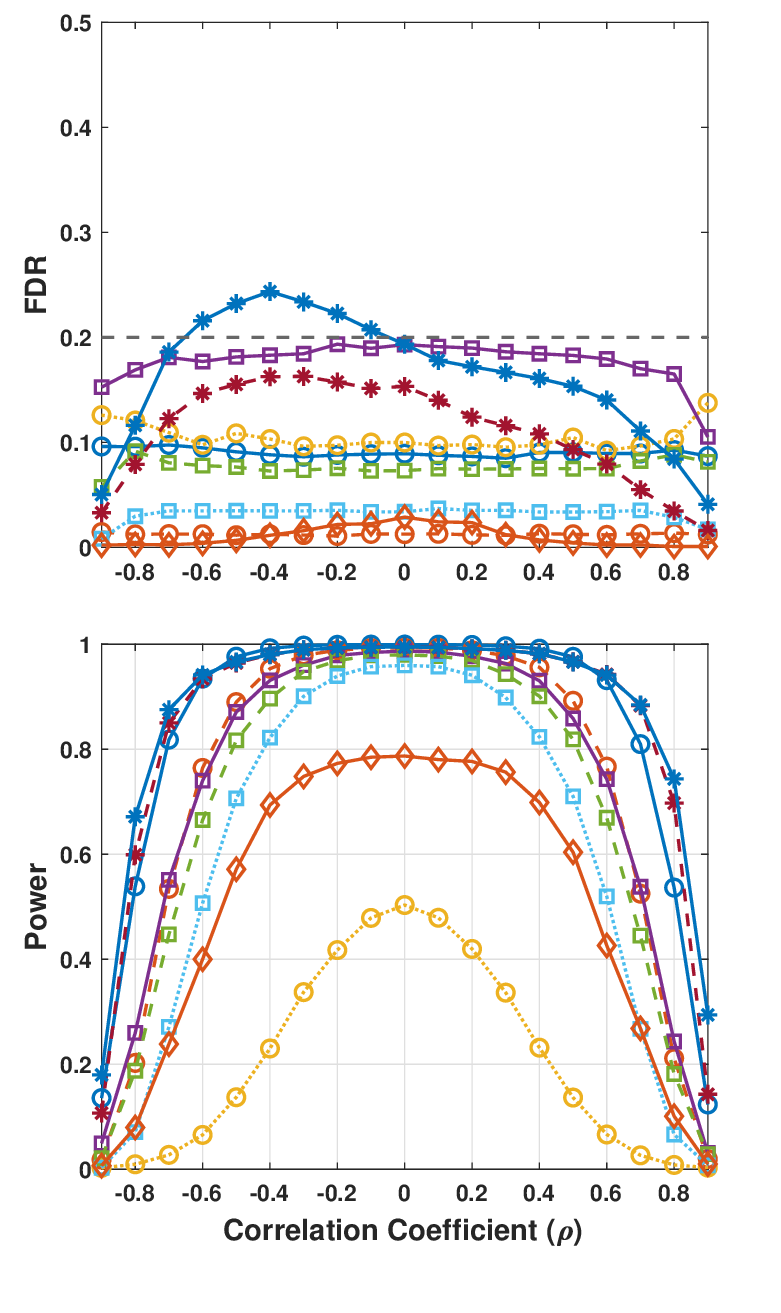}
\label{fig:subfig2}}
\caption{FDR and power vs Correlation coefficient $\rho$. In this experiment, the coefficients corresponding to the false nulls are set to be $\beta_a = 8$.
}
\label{fig:rho20}
\end{figure}

 The plots are based on averaging 200 trials and the power is defined as follows,
\begin{equation}
     \text{Power}:=\frac{1}{k}\,\mathbb{E}\big|\hat\Rc\cap{\Breve\Hc_0}^c \big|\ .
\end{equation}
{Our simulations show that the FRPP knockoff procedure with $\eps=0.8$ outperforms the knockoff-assisted BH procedure, the S-OLS methods perform similarly to the BY procedure, and our heuristic methods outperform the original BH procedure which serves as a heuristic baseline. However, S-LASSO2 shows slight FDR violation in simulation (b) of Figure~\ref{fig:rho20}.}

\begin{figure}[h!]
\centering
\subfloat[][]{
\includegraphics[width=0.4\textwidth]{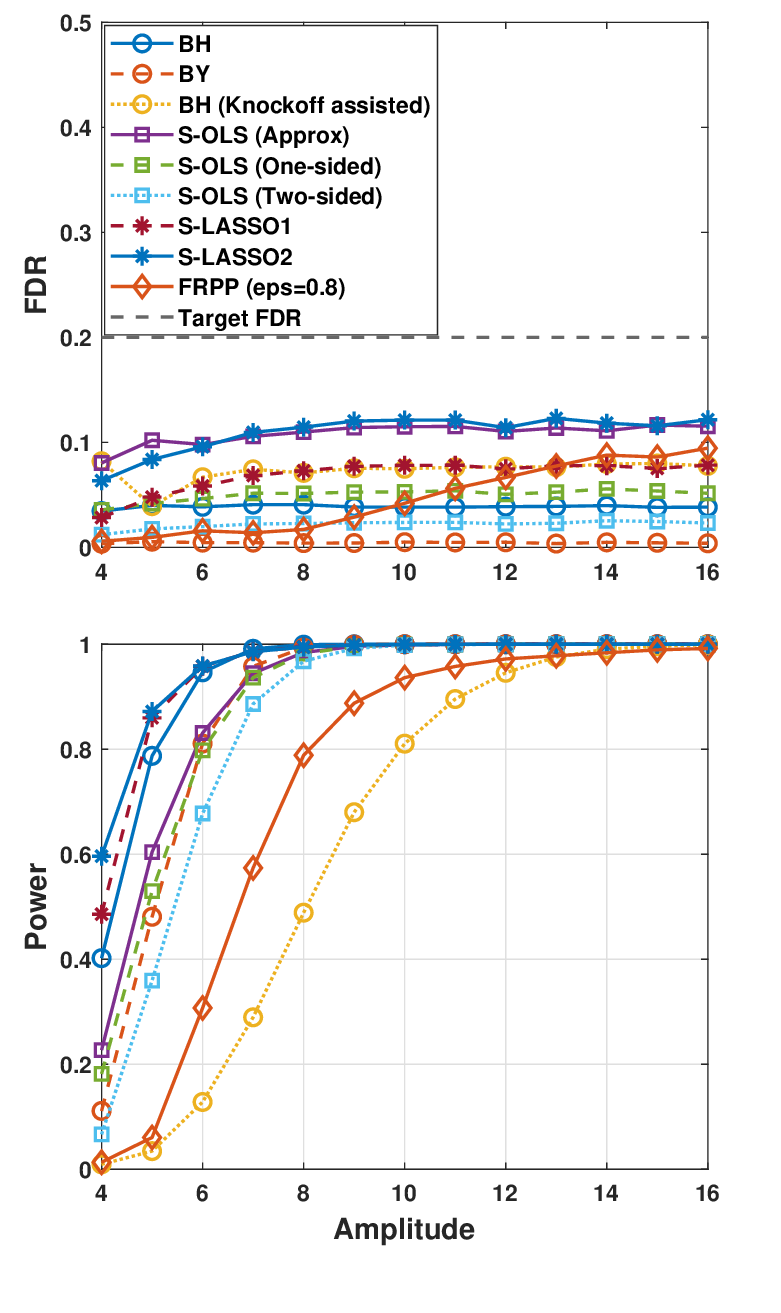}
\label{fig:subfig1}}
\subfloat[][]{
\includegraphics[width=0.4\textwidth]{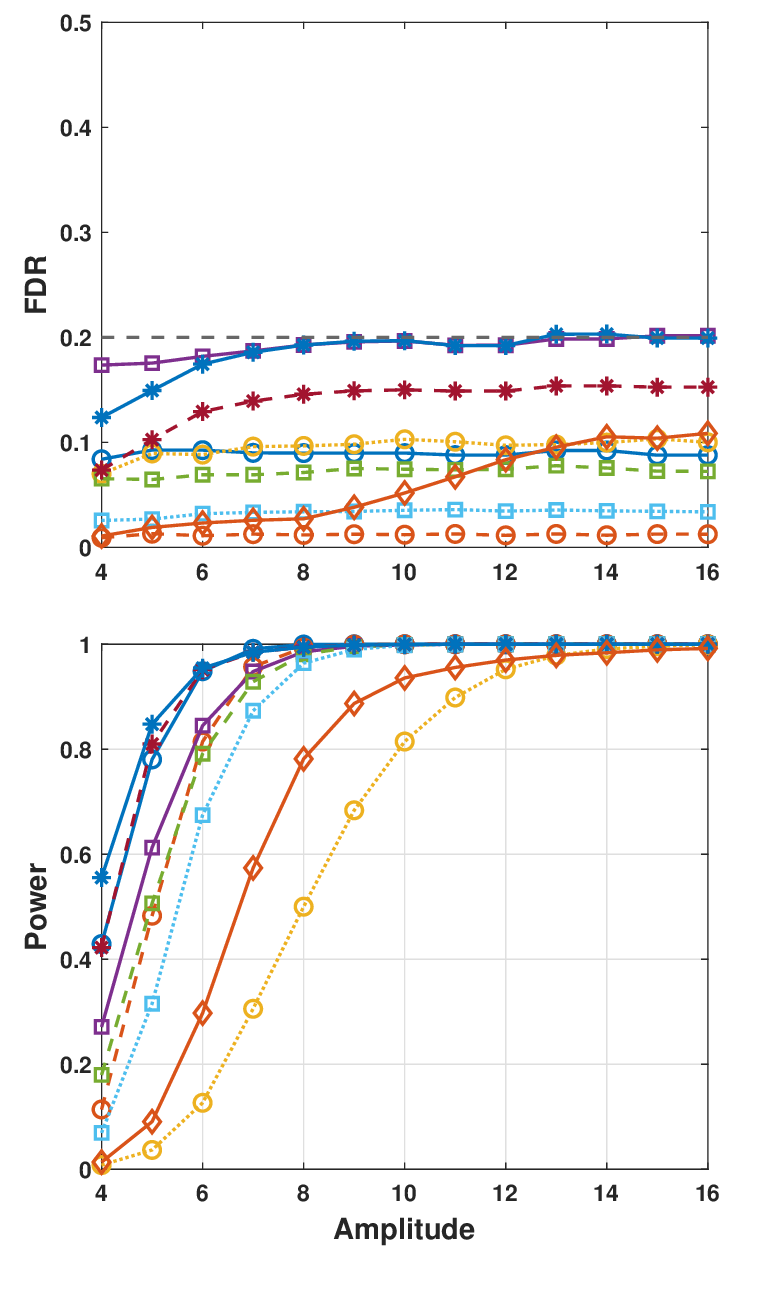}
\label{fig:subfig2}}
\caption{FDR and power vs the amplitude of the alternative coefficients. In this experiment the correlation coefficient is set to be $\rho = 0$ and all of the false nulls have the same underlying coefficient. 
}
\label{fig:amp20}
\end{figure}

\begin{remark}
In \cite{barber2015controlling}, it is suggested to set $s_i=\min(\\2\lambda_{\min}(\Sigmav),1)$ for equicorrelated knockoffs. However, this is not possible in the case of using OLS estimator as it would result in a singular {Gram matrix $\Gv$} when $2\lambda_{\min}(\Sigmav)\leq 1$ (See {the proof of} Lemma \ref{invlem}). Regarding the FRPP method, note that unlike the deterministic methods larger values of $s_i$ do not result in higher detection power necessarily, because the variance of the additive Laplacian noise is proportional to $s_i^2$. In our experiments, it turns out that for $\rho = 0$ (no correlation) case, the procedure reaches its highest power when $s_i=\min(\lambda_{\min}(\Sigmav),1)$. {However this is not the case for high correlations and $\lambda_{\min}(\Sigmav)$ would be too small. 
}
\end{remark}

\section{Discussion}
Fixed-X knockoff procedure is an elegant method for selective inference in linear models with FDR control guarantee. However, the composite extension of this method has not been developed yet. In this paper, we have investigated the fixed-X knockoff filter approach to the variable selection problem with composite nulls and under arbitrary dependencies among statistics. The knockoff inference procedure handles the dependencies between variables very naturally by computing model-based statistics. We have shown that this structure is still useful under the composite nulls and allows for the development of methods with theoretical FDR control guarantee. We have derived a full stochastic generalization of the knockoff procedure {by adding Laplace noise to the feature-response products and the method shows reasonable statistical power in simulations. Optimizing the noise variance and quantifying the relation between composite testing and the differentially private variable selection with knockoffs developed in \cite{pournaderi2021differentially} are left for future research.} We have shown that if we restrict ourselves to the ordinary least-squares estimates, the intuitive (and deterministic) method of shifting the estimates is theoretically valid. We have also derived a general bound on the FDR for cases where the original knockoff procedure is applied to a composite problem, without any additional assumptions.  

\section*{Acknowledgments}
\label{sec:ack}
The authors would like to thank the anonymous reviewers for their constructive comments and the Associate Editor for handling the submission.

\appendix

\section{Technical Lemmas}
\label{lemmas}
\begin{lemma} \label{invlem} 
If $\underset{1\leq i\leq p}{\max}\,s_i<2\lambda_{\min}(\Sigmav)$, then $\Gv = [\Xv\  \mathbf{\Tilde{X}}]^\top[\Xv\  \mathbf{\Tilde{X}}]$ is invertible.
\end{lemma}
\begin{proof}
From \eqref{eq:data_gram} recall that
\begin{equation}
    \Gv = 
    \begin{bmatrix} \boldsymbol{\Sigmav} &  \boldsymbol{\Sigmav}-\diag(\sv)\\\boldsymbol{\Sigmav}-\diag(\sv) &  \boldsymbol{\Sigmav}\end{bmatrix}\ ,\label{struct}
\end{equation}
where $\sv\in\Real_+^p$. 
We note,
\begin{equation*}
    \begin{bmatrix}
    {\Iv\quad-\Iv}\\{\Ov\qquad\Iv}
    \end{bmatrix}
    \begin{bmatrix}
    {\Av\quad \Bv}\\{\Bv\quad\Av}
    \end{bmatrix}
    \begin{bmatrix}
    {\Iv\quad\Iv}\\{\Ov\quad\Iv}
    \end{bmatrix} = 
    \begin{bmatrix}
    \Av-\Bv & \Ov\\\Bv & \Av+\Bv
    \end{bmatrix}\ .
\end{equation*}
Therefore, $\det(\Gv-\lambda \Iv) = \det\big(\diag(\sv)-\lambda \Iv\big)\det\big(2\Sigmav-\diag(\sv)-\lambda \Iv\big)$ and as a result, the set of eigenvalues of $\Gv$ is the union of the eigenvalues of $\diag(\sv)$ and $2\Sigmav-\diag(\sv)$. If $\underset{1\leq i\leq p}{\max}\,s_i<2\lambda_{\min}(\Sigmav)$ holds, then $\lambda_{\min}(2\Sigmav-\diag(\sv))>0$ which implies that $\Gv$ is positive definite and therefore, invertible.
\end{proof}

\begin{lemma} \label{inv}
If {$\underset{1\leq i\leq p}{\max}\,s_i<2\lambda_{\min}(\Sigmav)$}, then $\Gv^{-1}$ has the following structure.

\begin{equation*}
    \Gv^{-1} = 
    \begin{bmatrix} \mathbf{A} & \mathbf{A}-D^{-1}\\\mathbf{A}-D^{-1} & \mathbf{A}\end{bmatrix}\ ,
\end{equation*}
{where $\mathbf{A}=(2D-D\Sigmav^{-1}D)^{-1}$ and $D=\diag(\sv)$.}
\end{lemma}

\begin{proof}
From \eqref{eq:data_gram} recall that
\begin{equation}
    \Gv = 
    \begin{bmatrix} \boldsymbol{\Sigmav} & \boldsymbol{\Sigmav}-D\\\boldsymbol{\Sigmav}-D & \boldsymbol{\Sigmav}\end{bmatrix}\ .\label{struct}
\end{equation}
{Since $\underset{1\leq i\leq p}{\max}\,s_i<2\lambda_{\min}(\Sigmav)$, we get that $2D-D\Sigmav^{-1}D$ is invertible.}
From the inverse of $2\times 2$ block matrices we have
\begin{equation*}
    \Gv^{-1} = 
    \begin{bmatrix} \Av & -\Av(\Iv_p-D\Sigmav^{-1})\\-\Av(\Iv_p-D\Sigmav^{-1}) & \Av\end{bmatrix}\ .\nonumber
\end{equation*}
where $\Av=(2D-D\Sigmav^{-1}D)^{-1}$.
We observe,
\begin{align}
    \Av-\big(-\Av(\Iv_p-D\Sigmav^{-1})\big)={D}^{\,-1} \ ,
\end{align}
completing the proof.
\end{proof}

\begin{lemma}
Let $\begin{pmatrix}\boldsymbol{\gamma}_{p\times 1}\\\boldsymbol{\gamma'}_{p\times 1}\end{pmatrix}=[\Xv\ \mathbf{\Tilde{X}}]^\top\yv$. It holds that {$\big|\mathbb{E}(\gamma_j-\gamma_j')\big|\leq s_j\delta_j$} for all $j\in\Breve\Hc_0$.  \label{frp_mean}
\end{lemma}
\begin{proof}
We observe $\mathbb{E}\begin{pmatrix}\boldsymbol{\gamma}\\\boldsymbol{\gamma'}\end{pmatrix}= \Gv\begin{pmatrix}\betav\\\mathbf{0}_{p\times 1}\end{pmatrix}$. According to the structure of $\Gv$ \eqref{struct}, we get
\begin{subequations}
\begin{align*}
    \big|\mathbb{E}(\gamma_j-\gamma_j')\big| &= \bigg|\big(\Gv_{(j)}-\Gv_{(j+p)}\big)\begin{pmatrix}\betav\\ 0\end{pmatrix}\bigg| \\
    &= \bigg|\Big(\Sigmav_{(j)} - \big(\Sigmav_{(j)} - D_{(j)}\big)\Big)\betav\bigg| \\
    &=s_j|\beta_j|\leq s_j\delta_j\ ,
\end{align*}
\end{subequations}
where the index $(j)$ denotes the $j$-th row of the matrices and the inequality holds according to the definition of $\Breve{\mathsf{H}}_{0,j}$.
\end{proof}

\section{Proof of Theorem \ref{shifted_OLS}}
\begin{lemma} \label{thmols}
If we use estimates with distribution $\begin{pmatrix}\boldsymbol{\hat\beta}\\\boldsymbol{\hat\beta'}\end{pmatrix}\sim \mathcal{N}\big(\begin{pmatrix}\boldsymbol{\breve\beta}\\\boldsymbol{\breve\beta'}\end{pmatrix},\sigma^2\,\Gv^{-1}\big)$ to compute the statistics, then the following properties hold for any $1\leq j\leq p$. 
\smallskip

\noindent
\textbf{(I)} If $W_j$ depends on the estimator only through ${\hat\beta}_j$ and ${\hat\beta}_j'$, and  {$\mathsf{sgn}(W_j) = \mathsf{sgn}\big({\hat\beta}_j-{\hat\beta}_j'\big)$ or} $\mathsf{sgn}(W_j) = \mathsf{sgn}\big(\big|{\hat\beta}_j\big|-\big|{\hat\beta}_j'\big|\big)$, then
\begin{equation}
    \mathsf{sgn}(W_j)\indep \boldsymbol{W}_{-j}\Big| \big\{\hat\beta_{j},\hat\beta'_{j}\big\} , \label{OLS_indep_2}
\end{equation}
where $\{\cdot,\cdot\}$ denotes an unordered pair.
\smallskip

\noindent
{\textbf{(II)} Recall the block structure of $\Gv^{-1}$ from Lemma~\ref{inv}. We assert that
\begin{equation} 
    \begin{pmatrix}\boldsymbol{\hat\beta}+\boldsymbol{\hat\beta'}\\\boldsymbol{\hat\beta}-\boldsymbol{\hat\beta'}\end{pmatrix}\sim\mathcal{N}\left(\begin{pmatrix}\boldsymbol{\breve\beta}+\boldsymbol{\breve\beta'}\\\boldsymbol{\breve\beta}-\boldsymbol{\breve\beta'}\end{pmatrix},\sigma^2\,\begin{bmatrix} 4\Av-2D^{-1} & \Ov\\\Ov & 2D^{-1} \end{bmatrix}\right). \label{beta_dist}
\end{equation}
with $D=\diag(\sv)$.}
\end{lemma}
\begin{proof}
\textbf{(I)} We compute the conditional distribution, 
\begin{equation*}
    \Big(\boldsymbol{\hat\beta}_{-j},\boldsymbol{\hat\beta}'_{-j}\Big|\big(\hat\beta_{j},\hat\beta'_{j})=({a},{b}\big)\Big)\sim\mathcal{N}\left(\boldsymbol{\xi}^{(-j)},\mathbf{C}^{(-j)}\right)\ .
\end{equation*}
Let $\mathbf{K}=\sigma^2\,\Gv^{-1}$. {Using the structure of $\Gv^{-1}$ (discussed in Lemma \ref{inv}), standard calculations reveal,}
\begin{equation*}
    \boldsymbol{\xi}^{(-j)} = \begin{pmatrix}\boldsymbol{\breve\beta}_{-j}\\\boldsymbol{\breve\beta}'_{-j}\end{pmatrix}+\begin{bmatrix}
{\boldsymbol{c}\quad\boldsymbol{c}}\\{\boldsymbol{c}\quad\boldsymbol{c}}
\end{bmatrix}\begin{pmatrix}a-\breve\beta_{j}\\b-\breve\beta'_{j}\end{pmatrix}\nonumber\ ,
\end{equation*}
\begin{equation*}
   \mathbf{C}^{(-j)} = \mathbf{K}_{-j}-\begin{bmatrix}
{Q\quad Q}\\{Q \quad Q}
\end{bmatrix}\ ,
\end{equation*}
with some $\boldsymbol{c}\in\Real^{^{p-1}}$ and symmetric $Q\in \Real^{(p-1)\times(p-1)}$ that does not depend on $a$ and $b$. We note that the conditional distribution depends on the pair $(a,b)$ only through their sum $a+b$, so it is free of the order of the pair. Thus, we get
\begin{equation*}
    \Big(\boldsymbol{\hat\beta}_{-j},\boldsymbol{\hat\beta}'_{-j}\Big|\hat\beta_{j},\hat\beta'_{j}\Big)\overset{d}{=}\Big(\boldsymbol{\hat\beta}_{-j},\boldsymbol{\hat\beta}'_{-j}\Big|\big\{\hat\beta_{j},\hat\beta'_{j}\big\}\Big) \ .
\end{equation*}
Therefore, {if  $\mathsf{sgn}(W_j) = \mathsf{sgn}\big({\hat\beta}_j-{\hat\beta}_j'\big)$ or $\mathsf{sgn}(W_j) = \mathsf{sgn}\big(\big|{\hat\beta}_j\big|-\big|{\hat\beta}_j'\big|\big)$, it holds that}
\begin{equation*}
    \Big(\boldsymbol{W}_{-j}\Big|\mathsf{sgn}(W_j),\big\{\hat\beta_{j},\hat\beta'_{j}\big\}\Big)\overset{d}{=}\Big(\boldsymbol{W}_{-j}\Big|\big\{\hat\beta_{j},\hat\beta'_{j}\big\}\Big) \ ,
\end{equation*}
which implies \eqref{OLS_indep_2} immediately.

\textbf{(II)} 
{The proof follows from Lemma~\ref{inv} and standard calculations.}
\end{proof}
\subsection{One-sided Test:}
\begin{proof}
Let $\Hc_0' = \{i:\beta_i\leq \delta_i\}$. According to \eqref{FDRGbound}, it is sufficient to show 
\begin{equation*}
    {\underset{\omega:W_j \neq 0}{\text{ess}\sup}}\ \frac{\mathbb{P}\Big\{W_j > 0 \ \Big|\ |W_j|, \boldsymbol{W}_{-j}\Big\}}{\mathbb{P}\Big\{W_j < 0 \ \Big|\ |W_j|, \boldsymbol{W}_{-j}\Big\}} \leq 1,\quad \text{all } j\in \Hc_0'\, .
\end{equation*}
{By hypothesis, $|W_j|$ depends on $\boldsymbol{\hat{\beta}}_\mathsf{S-OLS}$ only through the unordered pair $\{{\hat\beta}_j,{\hat\beta}_j'\}$. Therefore, using the tower property we get,
\begin{align*}
    \frac{\mathbb{P}\Big\{W_j > 0 \ \Big|\ |W_j|, \boldsymbol{W}_{-j}\Big\}}{\mathbb{P}\Big\{W_j < 0 \ \Big|\ |W_j|, \boldsymbol{W}_{-j}\Big\}}&=
    \frac{\mathbb{E}\bigg\{\mathbb{P}\Big(W_j > 0 \,\Big| \big\{\hat\beta_{j},\hat\beta'_{j}\big\}, \boldsymbol{W}_{-j}\Big)\bigg|\,|W_j|, \boldsymbol{W}_{-j}\bigg\}}{\mathbb{E}\bigg\{\mathbb{P}\Big(W_j < 0 \,\Big| \big\{\hat\beta_{j},\hat\beta'_{j}\big\}, \boldsymbol{W}_{-j}\Big)\bigg|\,|W_j|, \boldsymbol{W}_{-j}\bigg\}}
    \  a.s.
\end{align*}
where $\{\cdot,\cdot\}$ denotes an unordered pair. Hence, it will be sufficient to prove,
\begin{equation}
    \underset{\omega:W_j \neq 0}{\text{ess}\sup}\ \frac{\mathbb{P}\Big(W_j > 0\Big| \big\{\hat\beta_{j},\hat\beta'_{j}\big\},\boldsymbol{W}_{-j}\Big)}{\mathbb{P}\Big(W_j < 0 \,\Big|\big\{\hat\beta_{j},\hat\beta'_{j}\big\},\boldsymbol{W}_{-j}\Big)}\leq 1,\quad \text{all } j\in \Hc_0'\, ,\label{abstow}
\end{equation}
which reduces to showing 
\begin{equation}
     \underset{\omega:\hat\beta_{j} \neq \hat\beta'_{j}}{\text{ess}\sup}\ \frac{\mathbb{P}\Big(\hat\beta_{j} > \hat\beta'_{j} \,\Big| \big\{\hat\beta_{j},\hat\beta'_{j}\big\}\Big)}{\mathbb{P}\Big(\hat\beta_{j} < \hat\beta'_{j} \,\Big| \big\{\hat\beta_{j},\hat\beta'_{j}\big\}\Big)} \leq 1,\quad \text{all } j\in \Hc_0'\, ,\label{ineq:one-sided}
\end{equation}
according to $\mathsf{sgn}(W_j) = \mathsf{sgn}\big({\hat\beta}_j-{\hat\beta}_j'\big)$ and \eqref{OLS_indep_2}.
Now we note,
\begin{align*}
     \underset{\omega:\hat\beta_{j} \neq \hat\beta'_{j}}{\text{ess}\sup}&\ \frac{\mathbb{P}\Big(\hat\beta_{j} > \hat\beta'_{j} \,\Big| \big\{\hat\beta_{j},\hat\beta'_{j}\big\}\Big)}{\mathbb{P}\Big(\hat\beta_{j} < \hat\beta'_{j} \,\Big| \big\{\hat\beta_{j},\hat\beta'_{j}\big\}\Big)} = 
     \underset{\omega:\hat\beta_{j} \neq \hat\beta'_{j}}{\text{ess}\sup}\ \frac{\mathbb{P}\Big(\hat\beta_{j} - \hat\beta'_{j}>0 \,\Big| \hat\beta_{j}+\hat\beta'_{j},\big|\hat\beta_{j} - \hat\beta'_{j} \big|\Big)}{\mathbb{P}\Big(\hat\beta_{j} - \hat\beta'_{j}<0 \,\Big| \hat\beta_{j}+\hat\beta'_{j},\big|\hat\beta_{j} - \hat\beta'_{j} \big|\Big)}\\
     &\qquad\qquad\qquad = \underset{\omega:\hat\beta_{j} \neq \hat\beta'_{j}}{\text{ess}\sup}\ \frac{\mathbb{P}\Big(\hat\beta_{j} - \hat\beta'_{j}>0 \,\Big| \big|\hat\beta_{j} - \hat\beta'_{j} \big|\Big)}{\mathbb{P}\Big(\hat\beta_{j} - \hat\beta'_{j}<0 \,\Big| \big|\hat\beta_{j} - \hat\beta'_{j} \big|\Big)}\leq 1,\quad \text{all } j\in \Hc_0',
\end{align*}
where the last equality holds according to \eqref{beta_dist} and the inequality follows from $\hat\beta_{j}-\hat\beta'_{j}\,{\sim}\,\Norm\left(\beta_j-\delta_j,2\sigma^2/s_j\right)$ and the hypotheses $\beta_j\leq \delta_j$.}

\end{proof}
\subsection{{Approximate} Two-sided Test:}
\begin{proof}
According to Theorem \ref{bcthm}, it is sufficient to show
\begin{equation*}
    {\underset{\omega:W_j \neq 0}{\text{ess}\sup}}\ \frac{\mathbb{P}\Big\{W_j > 0, \, \hat\beta_j+\hat\beta_j' > 0\,\Big|\ |W_j|, \boldsymbol{W}_{-j}\Big\}}{\mathbb{P}\Big\{W_j < 0 \ \Big|\ |W_j|, \boldsymbol{W}_{-j}\Big\}} \leq 1,\quad \text{all } j\in \Breve\Hc_0\, ,
\end{equation*}
which reduces to showing 
\begin{equation*}
    \underset{\omega:W_j \neq 0}{\text{ess}\sup}\ \frac{\mathbb{P}\Big(W_j > 0, \,\hat\beta_j+\hat\beta_j' > 0\Big|\boldsymbol{W}_{-j}, \big\{\hat\beta_{j},\hat\beta'_{j}\big\}\Big)}{\mathbb{P}\Big(W_j < 0 \,\Big|\boldsymbol{W}_{-j}, \big\{\hat\beta_{j},\hat\beta'_{j}\big\}\Big)}\leq 1,\quad \text{all } j\in \Breve\Hc_0\, ,
\end{equation*}
{by the same argument that led to \eqref{abstow} in the proof of one-sided test.} 
Since $\hat\beta_{j}+\hat\beta'_{j}$ is a function of $\big\{\hat\beta_{j},\hat\beta'_{j}\big\}$, we get
\begin{align*}
     &\frac{\mathbb{P}\Big(W_j > 0, \,\hat\beta_j+\hat\beta_j' > 0\Big|\boldsymbol{W}_{-j}, \big\{\hat\beta_{j},\hat\beta'_{j}\big\}\Big)}{\mathbb{P}\Big(W_j < 0 \,\Big|\boldsymbol{W}_{-j}, \big\{\hat\beta_{j},\hat\beta'_{j}\big\}\Big)}\\&\qquad\qquad
     =\frac{\ind\Big\{\hat\beta_j+\hat\beta_j' > 0\Big\}\mathbb{P}\Big(W_j > 0\Big|\boldsymbol{W}_{-j}, \big\{\hat\beta_{j},\hat\beta'_{j}\big\}\Big)}{\mathbb{P}\Big(W_j < 0 \,\Big|\boldsymbol{W}_{-j}, \big\{\hat\beta_{j},\hat\beta'_{j}\big\}\Big)}\\&\qquad\qquad
     {=}\ind\Big\{\hat\beta_j+\hat\beta_j' > 0\Big\}\frac{\mathbb{P}\Big(W_j > 0\Big| \big\{\hat\beta_{j},\hat\beta'_{j}\big\}\Big)}{\mathbb{P}\Big(W_j < 0 \,\Big| \big\{\hat\beta_{j},\hat\beta'_{j}\big\}\Big)}\quad a.s.,
\end{align*}
where the last equality holds according to $\eqref{OLS_indep_2}$. Now we note,
\begin{align*}
      &\underset{\omega:W_j \neq 0}{\text{ess}\sup}\ \ind\Big\{\hat\beta_j+\hat\beta_j' > 0\Big\}\frac{\mathbb{P}\Big(W_j > 0\Big| \big\{\hat\beta_{j},\hat\beta'_{j}\big\}\Big)}{\mathbb{P}\Big(W_j < 0 \,\Big| \big\{\hat\beta_{j},\hat\beta'_{j}\big\}\Big)}\\
      &\overset{(*)}{=}\underset{\omega:W_j \neq 0}{\text{ess}\sup}\ \ind\Big\{\hat\beta_j+\hat\beta_j' > 0\Big\}\frac{\mathbb{P}\Big(\hat\beta_{j} > \hat\beta'_{j}\Big| \big\{\hat\beta_{j},\hat\beta'_{j}\big\}\Big)}{\mathbb{P}\Big(\hat\beta_{j} < \hat\beta'_{j} \,\Big| \big\{\hat\beta_{j},\hat\beta'_{j}\big\}\Big)}\\
     &\leq\underset{\omega:\hat\beta_{j} \neq \hat\beta'_{j}}{\text{ess}\sup}\ \frac{\mathbb{P}\Big(\hat\beta_{j} > \hat\beta'_{j}\Big| \big\{\hat\beta_{j},\hat\beta'_{j}\big\}\Big)}{\mathbb{P}\Big(\hat\beta_{j} < \hat\beta'_{j} \,\Big| \big\{\hat\beta_{j},\hat\beta'_{j}\big\}\Big)}\leq 1,\quad \text{all } j\in \Breve\Hc_0\, .
\end{align*}
 where $(\ast)$ holds since $\hat\beta_{j} > \hat\beta'_{j} \Leftrightarrow |\hat\beta_{j}| > |\hat\beta'_{j}|$ when $\hat\beta_j+\hat\beta_j' > 0$ and the last inequality holds according to \eqref{ineq:one-sided}. 
\end{proof}

\section{Proof of Theorem \ref{lap_mech}}
\begin{proof}
In this method, the Laplace noise $\Delta_{2p\times 1}$ is added to the feature-response products,
\begin{equation*}
    [\Xv\  \mathbf{\Tilde{X}}]^\top\yv+\Delta = \Gv\,\overline{\betav}+\Delta+ [\Xv\  \mathbf{\Tilde{X}}]^\top\wv\ ,
\end{equation*}
where {$\Delta_j{\sim}\text{Lap}\,(2\,s_{j\,(\mathsf{mod}\,p)}\,\delta_{j\,(\mathsf{mod}\,p)}/\eps)$, $1\leq j\leq 2p$, denote $2p$ independent random variables}, $s_i = 1-\Xv_i^\top\mathbf{\Tilde{X}}_i$, $\overline{\betav}=\begin{pmatrix}\betav\\\mathbf{0}_{p\times 1}\end{pmatrix}$, and $\mathbf{w}$ is the model noise. For simplicity in notation we re-write the equation as follows,
\begin{equation*}
    \begin{pmatrix}\boldsymbol{\kappa}_{p\times 1}\\\boldsymbol{\kappa}'_{p\times 1}\end{pmatrix} = \begin{pmatrix}\boldsymbol{\theta}\\\boldsymbol{\theta}'\end{pmatrix} + \begin{pmatrix}\boldsymbol{\alpha}\\\boldsymbol{\alpha}'\end{pmatrix}\ ,
\end{equation*}
where $\begin{pmatrix}\boldsymbol{\theta}\\\boldsymbol{\theta}'\end{pmatrix} = \Gv\overline{\betav}+ \Delta$ and $\begin{pmatrix}\boldsymbol{\alpha}\\\boldsymbol{\alpha}'\end{pmatrix}=[\Xv\  \mathbf{\Tilde{X}}]^\top\wv$. By Lemma~\ref{frp_mean}, we have $\big|\mathbb{E}(\theta_j-\theta_j')\big|\leq s_j\,\delta_j$ for all $j\in\Breve\Hc_0$. Therefore, $\big|\mathbb{E}(\theta_j-\theta_j')\big|+\big|\mathbb{E}(\theta_j'-\theta_j)\big|\leq 2\,s_j\,\delta_j$ and the following inequality holds for all $(x,y)\in\Real^2$ and $j\in\Breve\Hc_0$ according to the Laplace mechanism by \cite{dwork2006calibrating,dwork2014algorithmic}.
\begin{equation}
    f_{(\theta_j,\theta'_{j})}(x,y) \leq e^\eps\ f_{(\theta'_j,\theta_{j})}(x,y)\ , \label{lapmech}
\end{equation}
where $f_{(\theta_j,\theta'_{j})}(x,y)$ denotes the probability density function.

In order to show the FDR control, according to \eqref{FDRGbound}, it will be sufficient to prove that
\begin{equation*}
    \mathbb{P}\Big\{W_j >  0 \ \Big|\ |W_j|,\, \boldsymbol{W}_{-j}\Big\} \leq e^\eps\,\mathbb{P}\Big\{{W_j <  0} \ \Big|\ |W_j|,\, \boldsymbol{W}_{-j}\Big\}\quad a.s.,
\end{equation*}
for all $j\in\Breve\Hc_0$.
We note that $\big(|W_j|,\boldsymbol{W}_{-j}\big)$ is a function of $(\big\{\kappa_j,\kappa_j'\},\boldsymbol{\alpha}_{-j},\boldsymbol{\alpha}'_{-j},\\\boldsymbol{\theta}_{-j},\boldsymbol{\theta}'_{-j} \big)$, therefore, according to the tower property we get
\begin{multline*}
    \mathbb{P}\Big\{W_j > 0 \ \Big|\ |W_j|\,,\, \boldsymbol{W}_{-j}\Big\}=\\
    \mathbb{E}\bigg\{\mathbb{P}\Big(W_j>0\Big|\{\kappa_j,\kappa_j'\},\boldsymbol{\alpha}_{-j},\boldsymbol{\alpha}'_{-j},\boldsymbol{\theta}_{-j},\boldsymbol{\theta}'_{-j}\Big) \bigg|\ |W_j|,\boldsymbol{W}_{-j}\bigg\}\quad a.s.,
    \end{multline*}
where $\{\cdot,\cdot\}$ denotes an unordered pair. 
Therefore, it is sufficient to show
\begin{align*}
    &\mathbb{P}\Big({W_j>0}\Big|\{\kappa_j,\kappa_j'\},\boldsymbol{\alpha}_{-j},\boldsymbol{\alpha}'_{-j},\boldsymbol{\theta}_{-j},\boldsymbol{\theta}'_{-j}\Big) \leq\nonumber\\
    &\qquad\qquad e^\eps\, \mathbb{P}\Big({W_j<0}\Big|\{\kappa_j,\kappa_j'\},\boldsymbol{\alpha}_{-j},\boldsymbol{\alpha}'_{-j},\boldsymbol{\theta}_{-j},\boldsymbol{\theta}'_{-j}\Big)\quad a.s.\ .\nonumber
\end{align*}
 {We also note that $\mathsf{sgn}(W_j)$ is conditionally independent of $(\boldsymbol{\theta}_{-j},\boldsymbol{\theta}'_{-j})$. Hence, we only need to show 
\begin{equation*}
    \underset{\omega:W_j \neq 0}{\text{ess}\sup}\ \frac{\mathbb{P}\Big(W_j>0\Big|\{\kappa_j,\kappa_j'\},\boldsymbol{\alpha}_{-j},\boldsymbol{\alpha}'_{-j}\Big)}{\mathbb{P}\Big(W_j<0\Big|\{\kappa_j,\kappa_j'\},\boldsymbol{\alpha}_{-j},\boldsymbol{\alpha}'_{-j}\Big)} \leq e^\eps\, .
\end{equation*}
We note
\begin{align*}
    \frac{\mathbb{P}\Big(W_j>0\Big|\{\kappa_j,\kappa_j'\}=\{u,v\},\boldsymbol{\alpha}_{-j},\boldsymbol{\alpha}'_{-j}\Big)}{\mathbb{P}\Big(W_j<0\Big|\{\kappa_j,\kappa_j'\}=\{u,v\},\boldsymbol{\alpha}_{-j},\boldsymbol{\alpha}'_{-j}\Big)}&\leq\max\Bigg\{\frac{\frac{g(u,v)}{g(u,v)+g(v,u)}}{\frac{g(v,u)}{g(u,v)+g(v,u)}},\frac{\frac{g(v,u)}{g(u,v)+g(v,u)}}{\frac{g(u,v)}{g(u,v)+g(v,u)}}\Bigg\}\nonumber
     \\& = \max\bigg\{\frac{g(u,v)}{g(v,u)},\frac{g(v,u)}{g(u,v)}\bigg\},
\end{align*}
almost everywhere on $\{(u,v,\omega):W_j\neq 0\}$,} with $g(\cdot,\cdot)$ denoting the shorthand for the conditional probability density function $g_{(\kappa_j,\kappa'_j)\big|(\boldsymbol{\alpha}_{-j},\boldsymbol{\alpha}'_{-j})}$. Now we note,
\begin{align*}
    &\nonumber{g_{(\kappa_j,\kappa'_j)\big|(\boldsymbol{\alpha}_{-j},\boldsymbol{\alpha}'_{-j})}{(u,v)}}
    ={g_{(\theta_j+\alpha_j,\theta'_{j}+\alpha'_{j})\big|(\boldsymbol{\alpha}_{-j},\boldsymbol{\alpha}'_{-j})}(u,v)}\\
    &\nonumber=\mathbb{E}_{(\theta_j,\theta'_{j})\big|(\boldsymbol{\alpha}_{-j},\boldsymbol{\alpha}'_{-j})}\Big\{ g_{(\alpha_j,\alpha'_{j})\big|(\boldsymbol{\alpha}_{-j},\boldsymbol{\alpha}'_{-j},\theta_{j},\theta'_{j})}(u-\theta_{j},v-\theta'_{j})\Big\}\\
    &\nonumber\overset{(*)}{=}{g_{(\theta_j+\alpha'_j,\theta'_{j}+\alpha_{j})\big|(\boldsymbol{\alpha}_{-j},\boldsymbol{\alpha}'_{-j})}(u,v)}\\
    &\nonumber={\mathbb{E}_{(\alpha_j,\alpha'_{j})\big|(\boldsymbol{\alpha}_{-j},\boldsymbol{\alpha}'_{-j})}\Big\{g_{(\theta_j,\theta'_{j})\big|(\boldsymbol{\alpha},\boldsymbol{\alpha}')}(u-\alpha'_{j},v-\alpha_{j})\Big\}}\\
     &\nonumber\overset{(**)}{\leq} {e^\eps}{g_{(\theta'_j+\alpha'_j,\theta_{j}+\alpha_{j})\big|(\boldsymbol{\alpha}_{-j},\boldsymbol{\alpha}'_{-j})}(u,v)}
    ={e^\eps}{g_{(\kappa_j,\kappa'_j)\big|(\boldsymbol{\alpha}_{-j},\boldsymbol{\alpha}'_{-j})}{(v,u)}},
\end{align*}
where 
$(\ast)$ holds according to the independence of $(\theta_j,\theta'_j)$ and $(\boldsymbol{\alpha},\boldsymbol{\alpha}')$, and the fact that
{
\begin{equation*}
    \Big(\alpha'_j,\boldsymbol{\alpha}_{-j}^\top,\alpha_j,{\boldsymbol{\alpha}'}^\top_{-j}\Big)\overset{d}{=} \Big(\alpha_j,\boldsymbol{\alpha}_{-j}^\top,\alpha'_j,{\boldsymbol{\alpha}'}^\top_{-j}\Big)
\end{equation*}
}since $\begin{pmatrix}\boldsymbol{\alpha}\\\boldsymbol{\alpha}'\end{pmatrix}=[\Xv\  \mathbf{\Tilde{X}}]^\top\wv\sim\mathcal{N}(0,{\sigma^2\Gv})$. The inequality marked with $(\ast\ast)$ follows from the independence of $(\theta_j,\theta'_j)$ and $(\boldsymbol{\alpha},\boldsymbol{\alpha}')$ and \eqref{lapmech}. {Therefore,
\begin{equation*}
    \underset{u,v,\omega}{\text{ess}\sup}\ \max\bigg\{\frac{g(u,v)}{g(v,u)},\frac{g(v,u)}{g(u,v)}\bigg\} \leq e^\eps,
\end{equation*}
completing the proof.}
\end{proof}
\section{Proof of Theorem \ref{boundthm}}
{The following lemma is equivalent to Lemma~1 in the supplementary material of \cite{barber2019knockoff}. We provide our alternative proofs, which follow techniques that are in line with that of Lemma~\ref{thmols} and Theorem~\ref{lap_mech}, for the sake of completeness.}
\begin{lemma}
 Let $\begin{pmatrix}\boldsymbol{\gamma}_{p\times 1}\\\boldsymbol{\gamma'}_{p\times 1}\end{pmatrix}=[\Xv\ \mathbf{\Tilde{X}}]^\top\yv$. For fixed $\Xv$, any anti-symmetric $\hat{\boldsymbol{W}} = \mathcal{W}\left([\Xv\ \mathbf{\Tilde{X}}]^\top[\Xv\ \mathbf{\Tilde{X}}],[\Xv\ \mathbf{\Tilde{X}}]^\top\yv\right)$ statisfies the following properties. \label{frp_lemma}

\smallskip
\noindent
\textbf{(I)} For all $1\leq j\leq p$, we have $\big(\gamma_j-\gamma_j'\big)\indep \big(\boldsymbol{\gamma}_{-j},\boldsymbol{\gamma}'_{-j}\big)$ and
\begin{equation}
    \mathsf{sgn}(W_j)\indep\mathsf{sgn}(\boldsymbol{W}_{-j})\,\Big|\big\{\boldsymbol{\gamma},\boldsymbol{\gamma}'\big\}_{vec}\ ,\label{frpindep}
\end{equation}
where {$\{\cdot,\cdot\}_{vec}$} denotes a vector of unordered pairs and $\mathsf{sgn}(\cdot)$ operates coordiante-wise.

\smallskip

\noindent
\textbf{(II)} For all $j\in\Breve\Hc_0$,
\begin{multline*}
     \mathbb{P}\Big\{W_j > 0 \,\Big|\{\gamma_{j},\gamma_j'\},(\boldsymbol{\gamma}_{-j},\boldsymbol{\gamma}'_{-j})\Big\}\leq\\ \quad \exp\left({\frac{\delta_j}{\sigma^2}\big|\gamma_{j}-\gamma_{j}'\big|}\right){\mathbb{P}\Big\{W_j < 0 \,\Big|\{\gamma_{j},\gamma_j'\},(\boldsymbol{\gamma}_{-j},\boldsymbol{\gamma}'_{-j})\Big\}}\nonumber \ a.s.,
\end{multline*}
where $\{\cdot,\cdot\}$ denotes an unordered pair.
\end{lemma}

\begin{proof}
\textbf{(I)} 
We compute the conditional distribution $(\gamma_j,\gamma'_j)\big|(\boldsymbol{\gamma}_{-j},\boldsymbol{\gamma}'_{-j})$. We note that $\begin{pmatrix}\boldsymbol{\gamma}\\\boldsymbol{\gamma'}\end{pmatrix}\sim \mathcal{N}\big(\begin{pmatrix}\boldsymbol{\breve\gamma}\\\boldsymbol{\breve\gamma'}\end{pmatrix},\sigma^2\,\Gv\big)$ where $\begin{pmatrix}\boldsymbol{\breve\gamma}\\\boldsymbol{\breve\gamma'}\end{pmatrix}=\Gv\begin{pmatrix}\betav\\\mathbf{0}_{p\times 1}\end{pmatrix}$ according to \eqref{linear_model}. Therefore we have,
\begin{equation}
    \Big((\gamma_j,\gamma'_j)\big|(\boldsymbol{\gamma}_{-j},\boldsymbol{\gamma}'_{-j})=(\boldsymbol{a},\boldsymbol{b})\Big)\sim\mathcal{N}(\boldsymbol{\eta}^{(j)},\sigma^2\,\mathbf{R}^{(j)})\ ,\label{conddist}
\end{equation}
where standard calculations together with one application of Lemma~\ref{inv} reveal,
\begin{align}
    \boldsymbol{\eta}^{(j)} &= \begin{pmatrix}\eta_j\\\eta_j'\end{pmatrix}:=\begin{pmatrix}\breve\gamma_j\\\breve\gamma'_j\end{pmatrix}+\begin{bmatrix}
{\boldsymbol{c_j}\quad\boldsymbol{c_j}}\\{\boldsymbol{c_j}\quad\boldsymbol{c_j}}
\end{bmatrix}^\top\begin{pmatrix}\boldsymbol{a}-\boldsymbol{\breve\gamma}_{-j}\\\boldsymbol{b}-\boldsymbol{\breve\gamma}'_{-j}\end{pmatrix}\label{condmean}\\ 
&=\begin{pmatrix}\breve\gamma_j+\boldsymbol{c_j}^\top(\boldsymbol{a}+\boldsymbol{b}-\boldsymbol{\breve\gamma}_{-j}-\boldsymbol{\breve\gamma}'_{-j}) \\\breve\gamma'_j+\boldsymbol{c_j}^\top(\boldsymbol{a}+\boldsymbol{b}-\boldsymbol{\breve\gamma}_{-j}-\boldsymbol{\breve\gamma}'_{-j}) \end{pmatrix}
= \begin{pmatrix}\breve\gamma_j+d_j\\\breve\gamma'_j+d_j\end{pmatrix}\nonumber\ ,
\end{align}
\begin{equation}
    \mathbf{R}^{(j)} = \begin{bmatrix}\Xv_j \ \mathbf{\Tilde{X}}_j\end{bmatrix}^\top\begin{bmatrix}\Xv_j \ \mathbf{\Tilde{X}}_j\end{bmatrix}-\begin{bmatrix}
{e_j\quad e_j}\\{e_j\quad e_j} 
\end{bmatrix}\ , \label{condcov}
\end{equation}
with some $\boldsymbol{c}_j\in\Real^{^{p-1}}$ and $e_j\in\Real$ that does not depend on $\boldsymbol{a}$ or $\boldsymbol{b}$. 
{Elementary calculations using \eqref{condmean} and \eqref{condcov} yields
\begin{equation}
    \Big(\gamma_j-\gamma'_j\big|\boldsymbol{\gamma}_{-j},\boldsymbol{\gamma}'_{-j}\Big)\overset{d}{=}\big(\gamma_j-\gamma_j'\big)\ , \label{diffindep}
\end{equation}
establishing the first claim.} According to~\eqref{condmean} and \eqref{condcov}, the conditional distribution \eqref{conddist} only depends on $\boldsymbol{a}+\boldsymbol{b}$ (through $\boldsymbol{\eta}^{(j)}$). Hence, we have
\begin{equation*}
    \Big(\gamma_j,\gamma'_j\Big|\boldsymbol{\gamma}_{-j},\boldsymbol{\gamma}'_{-j}\Big)\overset{d}{=}\Big(\gamma_j,\gamma'_j\Big|\big\{\boldsymbol{\gamma}_{-j},\boldsymbol{\gamma}'_{-j}\big\}_{vec}\Big)\ ,\nonumber
\end{equation*}
We observe that the unordered pair $\{\gamma_{j},\gamma'_{j}\}$ is a function of $(\gamma_j,\gamma'_j)$. Therefore, by conditioning on it we get,
\begin{equation*}
    \Big(\gamma_j,\gamma'_j\,\Big|(\boldsymbol{\gamma}_{-j},\boldsymbol{\gamma}'_{-j}),\{\gamma_{j},\gamma'_{j}\}\Big)\overset{d}{=}\\\Big(\gamma_j,\gamma'_j\Big|\big\{\boldsymbol{\gamma},\boldsymbol{\gamma}'\big\}_{vec}\Big)\nonumber\ .
\end{equation*}
Therefore,
\begin{equation*}
    \Big(\mathsf{sgn}(W_{j})\,\Big|\mathsf{sgn}(\boldsymbol{W}_{-j}),\{\boldsymbol{\gamma},\boldsymbol{\gamma}'\big\}_{vec}\Big)\overset{d}{=}\Big(\mathsf{sgn}(W_{j})\Big|\big\{\boldsymbol{\gamma},\boldsymbol{\gamma}'\big\}_{vec}\Big)\nonumber\ ,
\end{equation*}
which implies \eqref{frpindep} immediately.

\textbf{(II)}
We note that the desired inequality holds {trivially for $\{\omega:W_j=0\}$. 
Also, $\mathbb{P}\left\{W_j < 0 \Big|\{\gamma_{j},\gamma_j'\},(\boldsymbol{\gamma}_{-j},\boldsymbol{\gamma}'_{-j})\right\}>0$ for $\{\omega:W_j\neq 0\}$.} Therefore, it is sufficient to prove
\begin{equation*}
     \frac{\mathbb{P}\Big\{W_j > 0 \,\Big|\{\gamma_{j},\gamma_j'\},(\boldsymbol{\gamma}_{-j},\boldsymbol{\gamma}'_{-j})\Big\}}{\mathbb{P}\Big\{W_j < 0 \,\Big|\{\gamma_{j},\gamma_j'\},(\boldsymbol{\gamma}_{-j},\boldsymbol{\gamma}'_{-j})\Big\}} \leq \exp\Big\{{\frac{\delta_j}{\sigma^2}\big|\gamma_{j}-\gamma_{j}'\big|}\Big\},
\end{equation*}
 {almost everywhere on $\{\omega:W_j \neq 0\}$} and all $j\in \Breve\Hc_0$. Let $h(\cdot,\cdot)$ denote $h_{(\gamma_j,\gamma'_j)\big|(\boldsymbol{\gamma}_{-j},\boldsymbol{\gamma}'_{-j})}\overset{\text{pdf}}{\sim} \mathcal{N}(\boldsymbol{\eta}^{(j)},\sigma^2\,\mathbf{R}^{(j)})$. Notice that on $\{\omega:W_j\neq 0\}$, swapping $\gamma_j$ and $\gamma'_j$ would switch the sign of $W_j$. Hence, we get
\begin{align}
     \frac{\mathbb{P}\Big\{W_j > 0 \,\Big|\{\gamma_{j},\gamma_j'\}=\{u,v\},(\boldsymbol{\gamma}_{-j},\boldsymbol{\gamma}'_{-j})\Big\}}{\mathbb{P}\Big\{W_j < 0 \,\Big|\{\gamma_{j},\gamma_j'\}=\{u,v\},(\boldsymbol{\gamma}_{-j},\boldsymbol{\gamma}'_{-j})\Big\}} &\leq\max\Bigg\{\frac{\frac{h(u,v)}{h(u,v)+h(v,u)}}{\frac{h(v,u)}{h(u,v)+h(v,u)}},\frac{\frac{h(v,u)}{h(u,v)+h(v,u)}}{\frac{h(u,v)}{h(u,v)+h(v,u)}}\Bigg\}\nonumber
     \\&= \max\bigg\{\frac{h(u,v)}{h(v,u)},\frac{h(v,u)}{h(u,v)}\bigg\} , \label{frpbound}
\end{align}
{almost everywhere on $\{(u,v,\omega):W_j\neq 0\}$}. {According to \eqref{condcov}, a ${\pi}/{4}$ clockwise rotation of the coordinate system implies}
\begin{align}
    \frac{h(u,v)}{h(v,u)} =\frac{h_1\big(\frac{1}{\sqrt{2}}(u-v)\big)h_2\big(\frac{1}{\sqrt{2}}(u+v)\big)}{h_1\big(\frac{1}{\sqrt{2}}(v-u)\big)h_2\big(\frac{1}{\sqrt{2}}(v+u)\big)}=\frac{h_1\big(\frac{1}{\sqrt{2}}(u-v)\big)}{h_1\big(\frac{-1}{\sqrt{2}}(u-v)\big)} \ ,\label{rat}
\end{align}
where {$h_1=h_{\frac{1}{\sqrt{2}}(\gamma_j-\gamma'_j)\big|(\boldsymbol{\gamma}_{-j},\boldsymbol{\gamma}'_{-j})}$ and $h_2=h_{\frac{1}{\sqrt{2}}(\gamma_j+\gamma'_j)\big|(\boldsymbol{\gamma}_{-j},\boldsymbol{\gamma}'_{-j})}$ denote probability density functions. By \eqref{diffindep}, we have $h_1\overset{\text{pdf}}{\sim}\Norm\big(\frac{1}{\sqrt{2}}\big(\breve\gamma_j-{\breve\gamma_j}'\big),\sigma^2_j \big)$} 
where $\sigma^2_j=s_j\,\sigma^2$ and $s_j = 1-\Xv_j^\top\mathbf{\Tilde{X}}_j$ since the columns of $\Xv$ are normalized. 
We now {continue from} \eqref{rat} to bound the RHS of \eqref{frpbound} under $\Breve{\mathsf{H}}_{0,j}$,
\begin{align*}
    \max\bigg\{\frac{h(u,v)}{h(v,u)},\frac{h(v,u)}{h(u,v)}\bigg\} &=\max\Bigg\{\frac{h_1\big(\frac{1}{\sqrt{2}}(u-v)\big)}{h_1\big(\frac{-1}{\sqrt{2}}(u-v)\big)},\frac{h_1\big(\frac{-1}{\sqrt{2}}(u-v)\big)}{h_1\big(\frac{1}{\sqrt{2}}(u-v)\big)}\Bigg\} \nonumber\\
    &\le\frac{\exp\bigg\{-\frac{1}{4\sigma_j^2}\Big(\big|\breve\gamma_j-\breve\gamma_j'\big|-\big|u-v\big|\Big)^2\bigg\}}{\exp\bigg\{-\frac{1}{4\sigma_j^2}\Big(\big|\breve\gamma_j-\breve\gamma_j'\big|+\big|u-v\big|\Big)^2\bigg\}}\nonumber\\
    &=\exp\bigg\{\frac{1}{\sigma_j^2}\big|\breve\gamma_j-\breve\gamma_j'\big|\big|u-v\big|\bigg\}\\
    &\leq\exp\bigg\{\frac{s_j\delta_j}{\sigma_j^2}\,\big|u-v\big|\bigg\}
    =\exp\bigg\{\frac{\delta_j}{\sigma^2}\big|u-v\big|\bigg\}\, ,
\end{align*}
where the second inequality is a consequence of $\big|\breve\gamma_j-\breve\gamma_j'\big|\leq s_j\delta_j$ under $\Breve{\mathsf{H}}_{0,j}$ (Lemma \ref{frp_mean}). 
\end{proof}

\begin{lemma}\label{ko_qual}
For all $j\in\Breve\Hc_0$ we have
\begin{equation*}
    {\mathbb{P}\Big\{W_j > 0\, ,\  \frac{\delta_j}{\sigma^2}\,\big|\gamma_j-\gamma_j'\big|\leq \eps \ \Big|\ |W_j|, \boldsymbol{W}_{-j}\Big\}}\leq e^\eps \, {\mathbb{P}\Big\{W_j < 0 \ \Big|\ |W_j|, \boldsymbol{W}_{-j}\Big\}},
\end{equation*}
almost surely.
\end{lemma}
\begin{proof}
\begin{align*}
    \mathbb{P}&\Big\{W_j > \, 0\,,\ \frac{\delta_j}{\sigma^2}\,\big|\gamma_j-\gamma_j'\big|\leq \eps \ \Big|\ |W_j|, \boldsymbol{W}_{-j}\Big\} \\= &\,\mathbb{E}\bigg\{\mathbb{P}\Big(W_j > 0\, ,\  \frac{\delta_j}{\sigma^2}\,\big|\gamma_j-\gamma_j'\big|\leq \eps \,\Big| \{\gamma_{j},\gamma'_{j}\},(\boldsymbol{\gamma}_{-j},\boldsymbol{\gamma}'_{-j})\Big)\bigg|\,|W_j|,\boldsymbol{W}_{-j}\bigg\}\nonumber\\=&\,\mathbb{E}\bigg\{\ind\Big\{\frac{\delta_j}{\sigma^2}\,\big|\gamma_j-\gamma_j'\big|\leq \eps \Big\}\mathbb{P}\Big(W_j > 0 \, \Big| \{\gamma_{j},\gamma'_{j}\},(\boldsymbol{\gamma}_{-j},\boldsymbol{\gamma}'_{-j})\Big)\bigg|\,|W_j|,\boldsymbol{W}_{-j}\bigg\}\nonumber\\\leq &\,\mathbb{E}\bigg\{\ind\Big\{\frac{\delta_j}{\sigma^2}\,\big|\gamma_j-\gamma_j'\big|\leq \eps \Big\}\exp\Big({\frac{\delta_j}{\sigma^2}\big|\gamma_j-\gamma_j'\big|}\Big)\mathbb{P}\Big(W_j < 0 \, \Big| \{\gamma_{j},\gamma'_{j}\},(\boldsymbol{\gamma}_{-j},\boldsymbol{\gamma}'_{-j})\Big)\nonumber\\
    &\qquad\qquad\qquad\qquad\qquad\qquad\qquad\qquad\qquad\qquad\qquad\qquad\qquad\bigg|\,|W_j|,\boldsymbol{W}_{-j}\bigg\}\nonumber\\\leq & \,e^{\eps}\,\mathbb{E}\bigg\{\mathbb{P}\Big(W_j < 0  \Big| \{\gamma_{j},\gamma'_{j}\},(\boldsymbol{\gamma}_{-j},\boldsymbol{\gamma}'_{-j})\Big)\bigg||W_j|,\boldsymbol{W}_{-j}\bigg\}\nonumber\\
    =&\,e^\eps\, \mathbb{P}\Big\{W_j < 0 \ \Big|\ |W_j|, \boldsymbol{W}_{-j}\Big\}\nonumber \quad a.s.,
\end{align*}
where $\{\cdot,\cdot\}$ denotes an unordered pair, the first inequality holds according to Lemma \ref{frp_lemma} (II), and we used the tower property since $\big(|W_j|,\boldsymbol{W}_{-j}\big)$ is a function of $\big(\{\gamma_{j},\gamma'_{j}\},(\boldsymbol{\gamma}_{-j},\boldsymbol{\gamma}'_{-j}) \big)$.
\end{proof}
Hence \eqref{naivebnd} holds according to Theorem \ref{bcthm}.

 

\bibliographystyle{elsarticle-num-names} 
\bibliography{cas-refs}





\end{document}